\documentclass[11pt,a4paper,oneside]{article}

\usepackage{amsmath}
\usepackage{amsfonts}
\usepackage{amssymb}
\usepackage{amsthm}
\usepackage{latexsym}
\usepackage[active]{srcltx}
\usepackage{mathrsfs}

\usepackage{enumerate}
\usepackage{color}

 \sloppypar

\newtheorem{theorem}{Theorem}[section]

\newtheorem{proposition}[theorem]{Proposition}

\newcommand{\im}{\text{Im }}

\newcommand{\tr}{\text{tr}}

\renewcommand{\H}{{\mathcal H}}

\def\ga{\gamma}
\def\GA{\Gamma}

\def\va{\varphi}

\def\g{\mathfrak{g}}
\def\a{\mathfrak{a}}

\def\ll{\mathfrak{l}}

\def\s{\mathfrak{s}}

\def\ga{\gamma}

\def\ps{\psi}

\def\ol#1{\overline{#1}}
\def\nn{\nonumber}

\def\R{{\mathbb R}}
\def\C{{\mathbb C}}
\def\N{{\mathbb N}}

\def\Z{{\mathbb Z}}

\def\B{{\mathcal B}}
\def\D{{\mathcal D}}
\def\F{{\mathcal F}}

\def\H{{\mathcal H}}
\def\P{{\mathcal P}}
\def\K{{\mathcal K}}

\def\U{{\mathcal U}}

\def\tr{{\text tr}}

\def\id{\text{id}}
\def\iy{\infty}

\def\wh{\widehat}

\def\ol#1{\overline{#1}}

\def\hb#1{\hbox{#1}}

\def\no#1{\Vert #1\Vert }

\def\ker#1{\hb{ker}(#1)}

\def\res#1{_{\vert #1}}

\def\hb #1{\hbox{#1}}

\def\hb#1{\hbox{#1}}

\def\ker#1{\hb{ker}(#1)}

\def\im#1{\hb{im}(#1)}

\def\noop#1{\Vert #1\Vert_{\rm op}}

\def\L1#1{L^1(#1)}

\def\L#1#2{L^{#1}(#2)}

\def\Im{\mathrm{\, Im \,}}
\def\Re{\mathrm{\, Re\, }}

\def\lef({\left(}
\def\rig){\right)}

 \sloppypar

\usepackage{a4wide}

\newtheorem{definition}[theorem]{Definition}

\newtheorem{lemma}[theorem]{Lemma}
\newtheorem{remark}[theorem]{Remark}
\newtheorem{convention}[theorem]{Convention}

\def\im{Im}

\newcommand{\quer}{\overline}

\def\tr{\mathrm{\, tr \,}}

\def\a{\mathfrak a}

\begin{document}

\title{The $C^*$-algebra of $SL(2,\R)$}
\author{Janne-Kathrin G\"unther}
\maketitle
~\\

\begin{abstract}
The $C^*$-algebra of the group $SL(2, \R)$ is charac\-terized using the operator valued Fourier transform. In 
particular, it is shown by explicit computations, that the Fourier transform 
of this $C^*$-algebra fulfills the norm controlled dual limit property.
\end{abstract}

\section{Introduction}
\label{intro}
In this article, the structure of the $C^*$-algebra of the group $G=SL(2, \R)$ will be analyzed. \\
The structure of the group $C^*$-algebras is already known for certain classes of Lie groups: The $C^*$-algebras of the Heisenberg and the thread-like Lie groups have been analyzed in [\ref{ludwig-turowska}] and the $C^*$-algebras of the $ax+b$-like groups in [\ref{lin-ludwig}]. Furthermore, the $C^*$-algebras of the 5-dimensional nilpotent Lie groups have been determined in [\ref{hedi}] and H.Regeiba characterized the $C^*$-algebras of all 6-dimensional nilpotent Lie groups in his doctoral thesis (see [\ref{hedidr}]). Just recently, the $C^*$-algebras of connected real two-step nilpotent Lie groups have been analyzed in [\ref{janne}]. \\
For semisimple Lie groups, there is no explicit description of their group $C^*$-algebras given in literature. However, for those semisimple Lie groups whose unitary dual is classified, the procedure of the determination of the group $C^*$-algebra used in this article might be successfully applied in a similar way. 
A characterization of reduced group $C^*$-algebras of semisimple Lie groups can be found in [\ref{was}].\\
In the present paper, the group $C^*$-algebra of $SL(2, \R)$ shall be described very explicitely. It will be shown that it is charac\-terized by some conditions which are called "norm controlled dual limit conditions" and which will be given in Section \ref{cond} below. In an abstract existence result in [\ref{bel-bel-lud}] these conditions are shown to hold true for all simply connected connected nilpotent Lie groups. They are explicitely checked for all 5- and 6-dimensional nilpotent Lie groups (see [\ref{hedi}]), for the Heisenberg Lie groups and the thread-like Lie groups (see [\ref{ludwig-turowska}]) and for the connected real two-step nilpotent Lie groups (see [\ref{janne}]). \\
At the beginning of this article, some notations and important facts which are needed in order to determine the $C^*$-algebra of $SL(2, \R)$ will be recalled. In Section \ref{cond} the above mentioned conditions which are characterizing a group $C^*$-algebra will be defined. The main result of this article, namely the compliance of the group $C^*$-algebra of $SL(2, \R)$ with these conditions, will be formulated and its proof will be accomplished in the following sections. Section \ref{unitary dual} is about the unitary dual of $SL(2, \R)$ and its topology and in Section \ref{cond 3a} and \ref{cond 3b}, the above specified conditions will be verified for the group $G=SL(2, \R)$. Finally, in Section \ref{result}, an alternative version of the result about the $C^*$-algebra of $SL(2, \R)$ will be presented and the concrete structure of $C^*(G)$ will be given.

\section{Preliminaries} \label{pre-sl2r}
\subsection{General definitions}\label{gen-def-sl2r}
\begin{definition}[Fourier transform] 
~\\
The Fourier transform $\F(a)=\hat{a}$ of an element $a$ of a $C^*$-algebra $C$ is defined in the following way: One chooses for every $\gamma \in \wh{C}$, the unitary dual of $C$, a representation $(\pi_{\gamma},\H_{\gamma})$ in the equivalence class of $\gamma$ and defines
$$\F(a)(\gamma):=\pi_{\gamma}(a) \in \B(\H_{\gamma}).$$
Then $\F(a)$ is contained in the algebra of all bounded operator fields over $\wh{C}$
$$l^{\infty}\big(\wh{C} \big)= \Big\{ \phi= \big( \phi(\pi_{\ga}) \in \B(\H_{\ga})\big)_{\ga \in \wh{C}}~|~ \| \phi \|_{\infty}:= \sup \limits_{\ga \in \wh{C}} \|\phi(\pi_{\ga})\|_{op}< \infty \Big\}$$
and the mapping $$\F: C \to l^{\infty}\big(\wh{C}\big),~ a \mapsto \hat{a}$$
is an isometric $*$-homomorphism. 
\end{definition}

\begin{definition}[Properly converging sequence]
~\\
A sequence $(t_k)_{k \in \N} $ in a topological space is called properly converging, if $(t_k)_{k \in \N}$ has limit points and if every subsequence of $(t_k)_{k \in \N}$ has the same limit set as $(t_k)_{k \in \N}$. 
\end{definition}

Recall that the $C^*$-algebra of a locally compact group $G$ is defined as the completion of the convolution algebra $L^1(G)$ with respect to the $C^*$-norm of $L^1(G)$, i.e. 
$$C^*(G):= \quer{L^1(G)}^{\| \cdot \|_{C^*(G)}}~~~\text{with}~~~\|f \|_{C^*(G)}:=\sup \limits_{\pi \in \wh{G}}\| \pi(f) \|_{op}.$$

A well-known result, that can be found in [\ref{dix}], states that the unitary dual of $C^*(G)$ coincides with the unitary dual of $G$: $$\wh{C^*(G)}=\wh{G}.$$

The unitary dual $\wh{G}$ has a natural topology which can be characterized in the following way: 

\begin{theorem}[Topology of the dual space] \label{top-dual}
~\\
Let $(\pi_k,\mathcal{H}_{\pi_k})_{k\in\mathbb{N}}$ be a family of irreducible unitary representations of a locally compact group $G$. Then $(\pi_k)_{k \in \N}$ converges to $(\pi, \H_{\pi})$ in $\widehat{G}$ if and only if for some non-zero (respectively for every) vector $\xi$ in $\mathcal{H}_{\pi}$, for every $k\in\N$ there exists $\xi_k\in\mathcal{H}_{\pi_k}$ such that the sequence of matrix coefficients $\big(\big\langle \pi_k(\cdot) \xi_k, \xi_k \big\rangle\big)_{k \in \N}$ converges uniformly on compacta to the matrix coefficient
$\big\langle \pi(\cdot) \xi, \xi \big\rangle$.
\end{theorem}
The proof of this theorem can be found in [\ref{dix}].

\subsection{The Lie group $SL(2, \R)$}\label{prel-sl2r}
From now on let
$$G:= SL(2, \R)=\{A \in M(2, \R)|~\det A=1\}$$
and let 
$$K:=SO(2)
=\bigg\{k_{\va}:=\begin{pmatrix}
\cos{\va} & -\sin{\va}  \\
\sin{\va} & \cos{\va}
\end{pmatrix} \bigg|~\va \in [0, 2 \pi) \bigg\}$$
be its maximal compact subgroup. Furthermore, define the one-dimensional subgroups
$$N:=\bigg\{\mu_x:=\begin{pmatrix}
1 & x  \\
0 & 1
\end{pmatrix} \bigg|~x \in \R \bigg\}
~~~\text{and}~~~
A:=\bigg\{a_t:=\begin{pmatrix}
e^{\frac{t}{2}} & 0  \\
0 & e^{-\frac{t}{2}}
\end{pmatrix} \bigg|~t \in \R \bigg\}$$
of $G$. Let
$$\g=\s \ll(2,\R)= \{A \in M(2, \R)|~\tr A=0\}$$
be the Lie algebra of $G$. \\
~\\
From the Iwasawa decomposition, $G=KAN$ and thus, for every $g \in G$ there exist $\kappa(g) \in K$, $\mu \in N$ and $H(g) \in \a$, where $\a=\bigg\{\begin{pmatrix}
t & 0  \\
0 & -t
\end{pmatrix} \Big|~t \in \R \bigg\}$ is the Lie algebra of $A$, such that
$$g=\kappa(g)e^{H(g)}\mu.$$
~\\
Moreover, define on $\a$ the mappings $\rho$ and $\nu_z$ for $z \in \C$ as
$$\rho \begin{pmatrix}
t & 0  \\
0 & -t
\end{pmatrix}:=t~~~\text{and}~~~
\nu_z \begin{pmatrix}
t & 0  \\
0 & -t
\end{pmatrix}:=zt~~~\forall t \in \R.$$

Furthermore, let
$$L^2(K)_+:=\big\{f \in L^2(K, \C)|~f(k)=f(-k)~\forall k \in K \big\}~~~\text{and}$$
$$L^2(K)_-:=\big\{f \in L^2(K, \C)|~f(k)=-f(-k)~\forall k \in K \big\}$$
and define for every $u \in \C$ the representations $\P^{+,u}$ on $\H_{\P^{+,u}}:=L^2(K)_+$ and $\P^{-,u}$ on $\H_{\P^{-,u}}:=L^2(K)_-$ as
$$\P^{\pm,u}(g)f(k):=e^{-(\nu_{u}+\rho)H(g^{-1}k)}f \Big(\kappa \big(g^{-1}k \big)\Big)~~~\forall g \in G~\forall f \in L^2(K)_{\pm}~ \forall k \in K.$$

\begin{remark}\label{irred}
~\\
The representation $\big(\P^{+,u},\H_{\P^{+,u}} \big)$ is irreducible if and only if $u \not\in 2 \Z+1$ and the representation $\big(\P^{-,u},\H_{\P^{-,u}} \big)$ is irreducible if and only if $u \not\in 2 \Z$. \\
Furthermore, $\big(\P^{+,u},\H_{\P^{+,u}} \big)$ and $\big(\P^{-,u},\H_{\P^{-,u}} \big)$ are unitary for $u \in i \R$.
\end{remark}

For the proof see [\ref{knapp}], Chapter 2. 

\begin{convention}
~\\
Throughout this paper, by $L^2(K)$ and $C^{\iy}(K)$ is meant $L^2(K, \C)$ and $C^{\iy}(K, \C)$, respectively. 
\end{convention}

\begin{lemma}\label{intabsch}
~\\
For every function $f \in L^2(K)_{\pm}$ and every $g \in G$
$$\int \limits_K e^{-2 \rho H(g^{-1}k)}\Big| f \Big(\kappa \big(g^{-1}k \big) \Big) \Big|^2 dk~=~ \|f\|_{L^2(K)}^2.$$
\end{lemma}

The proof can be found in [\ref{knapp}], Chapter VII.2.

\begin{definition}[$n$-th isotypic component]
~\\
For a representation $(\tilde{\pi},\H_{\tilde{\pi}})$ of $K$ define for every $n \in \Z$ the $n$-th isotypic component or $K$-type of $\tilde{\pi}$ as
$$\H_{\tilde{\pi}}(n):=\big\{v \in \H_{\tilde{\pi}} |~\tilde{\pi}(k_{\va})v=e^{in \va}v ~\forall \va \in [0,2 \pi) \big\}.$$
A representation $(\tilde{\pi},\H_{\tilde{\pi}})$ of $G$ is called even (respectively odd), if $\H_{{\tilde{\pi}}\res{K}}(n)=\{0\}$ for all \mbox{odd $n$} (respectively for all even $n$).
\end{definition}

Every irreducible unitary representation of $G$ is even or odd. \\
Furthermore, the algebraic direct sum 
$$ \bigoplus \limits_{n \in \Z} \H_{\tilde{\pi}}(n)$$
is dense in $\H_{\tilde{\pi}}$.

\begin{remark}\label{even-odd}
~\\
By the definition of the Hilbert spaces $L^2(K)_{\pm}$ of $\P^{\pm,u}$ for $u \in \C$, it is easy to verify that $\P^{+,u}$ is even for every $u \in \C$ and that $\P^{-,u}$ is odd for every $u \in \C$.
\end{remark}

By the definition of the $n$-th isotypic component, one can remark that 
$$\H_{\P^{+,u}}(n)= \C \cdot e^{-in \cdot}~\text{for all even}~n \in \Z~~~\text{and}~~~\H_{\P^{-,u}}(n)= \C \cdot e^{-in \cdot}~\text{for all odd}~n \in \Z.$$

\begin{definition}[$p_n$]
~\\
Denote for $n \in \Z$ by $b_n(f)$ the $n$-th Fourier coefficient of $f \in L^2(K)_{\pm}$ which is defined as
$$b_n(f):= \frac{1}{|K|} \int \limits_K f(k_{\va}) e^{-in \va}d k_{\va},$$
and let 
$$p_n(f):=b_{-n}(f) e^{-in \cdot}.$$
\end{definition}

One can easily show that for every $u \in \C$ and for every $n \in \Z$ the operator $p_n$ is the projection from $\H_{\P^{\pm,u}}=L^2(K)_{\pm}$ to the $n$-th isotypic component of the representation $\P^{\pm,u}$.

\section{$C^*$-algebras with norm controlled dual limits} \label{cond}
In this section, the definition of a $C^*$-algebra with "norm controlled dual limits", which was mentioned in the introduction, will be given.

~\\
\begin{definition}\label{NCDL}
~\\
A $C^*$-algebra $C$ is called $C^*$-algebra with \textbf{norm controlled dual limits} if it fulfills the following conditions:
\begin{enumerate}[{\bf{-}}]
\item {\bf{Condition 1}}: Stratification of the unitary dual: 
\begin{enumerate}
\item There is a finite increasing family $ S_0\subset S_1\subset\ldots\subset S_r=\wh{C} $ of closed
subsets of the unitary dual $\wh{C} $ of $ C$ in such a way that for $ i \in \{1,\cdots, r\} $ the subsets $\GA_0:=S_0$ and
$ \GA_i:=S_i\setminus S_{i-1}$ are Hausdorff in their relative topologies and such that $S_0$ consists of all the characters of $C$.
\item For every $ i\in \{0,\cdots, r\} $ there is a Hilbert space $ \H_i $ and for every $ \ga\in \GA_i $ there is a concrete realization $ (\pi_\ga,\H_i) $ of $ \ga $ on the Hilbert space $ \H_i $.
\end{enumerate}
\item {\bf{Condition 2}}: CCR $C^*$-algebra: \\
$C$ is a separable CCR (or liminal) $C^*$-algebra, i.e. a separable $C^*$-algebra such that the image of every irreducible representation
$ (\pi,\H)  $ of $ C $ is contained in the algebra of compact operators $ \K(\H) $ (which implies that the image equals $\K(\H)$).
\item {\bf{Condition 3}}: Changing of layers: \\
Let $a \in C$. 
\begin{enumerate}
\item  The mappings $ \ga \mapsto \F (a)(\ga) $ are norm continuous on the different sets $ \GA_i $.
\item   For any  $ i \in \{0,\cdots, r\} $  and for any
converging sequence contained in $ \GA_i $ with limit set outside $ \GA_i $ (thus in $S_{i-1}$), there is a properly converging subsequence $\ol\ga=(\ga_j)_{j\in\N} $, as well as a constant $ c>0 $ and for every $ j\in\N $ an involutive linear mapping \mbox{$ \tilde{\nu}_j=\tilde{\nu}_{\ol\ga,j}: CB(S_{i-1})\to \B(\H_i)$}, {which is  bounded by $ c\|\cdot\|_{S_{i-1}} $} (uniformly in $j$), such that
\begin{eqnarray*}
&\lim \limits_{j\to\iy }\big\|
\F (a)(\ga_j)-\tilde{\nu}_{j}\big(\F (a)\res{S_{i-1}}\big)\big\|_{op}=0. 
\end{eqnarray*}
Here $CB(S_{i-1})$ is the $ * $-algebra of all the uniformly
bounded fields of operators \linebreak
$ \big(\ps(\ga)\in \B(\H_l) \big)_{\ga \in \GA_l, l=0,\cdots, i-1}$, which are
operator norm continuous on the subsets $ \GA_l$ for every $  l\in\{0,\cdots, i-1\} $, provided with the infinity-norm 
\begin{eqnarray*}
\no{\ps}_{S_{i-1}}:=\sup_{\ga\in S_{i-1}}\noop{\ps(\ga)}. 
\end{eqnarray*}
 \end{enumerate}
 \end{enumerate}
\end{definition}

\begin{theorem}\label{NCDL - sl2r}
~\\
The $C^*$-algebra of $G=SL(2, \R)$ has norm controlled dual limits.
\end{theorem}

\begin{remark}
~\\
Throughout the rest of this article, this theorem will be proved. Concrete subsets $\Gamma_i$ and $S_i$ of $\wh{C^*(G)} = \wh{G} $ will be defined and in Section \ref{cond 3b}, the mappings $\big(\tilde{\nu}_{j} \big)_{j \in \N}$ will be constructed.  
\end{remark}

The norm controlled dual limit conditions completely characterize the structure of a group $C^*$-algebra in the following sense: Taking the number $r$, the Hilbert spaces $\H_i$, the sets $\Gamma_i$ and $S_i$ for $i \in \{0,...,r\}$ (see Section \ref{gammas} and Section \ref{rep} for their construction) and the mappings $\big(\tilde{\nu}_{j} \big)_{j \in \N}$ (see Section \ref{cond 3b} for their construction) required in the above definition, by [\ref{hedi}], Theorem 3.5, one gets the result below for the $C^*$-algebra of $G=SL(2, \R)$: 

\begin{theorem}\label{theo1 - sl2r}
~\\
The $C^*$-algebra $C^*(G)$ of $G=SL(2, \R)$ is isomorphic (under the Fourier transform) to the set of all operator fields $ \va $ defined over $ \wh G $ such that
\begin{enumerate}\label{}
\item $ \va(\ga)\in \K(\H_i) $ for every $i \in \{1,...,r\}$ and every $ \ga\in\GA_i$.
\item $\va \in l^{\infty}(\wh{G})$.
\item  The mappings $ \ga \mapsto \va(\ga) $ are norm continuous on the different sets $ \GA_i $.
\item  For any sequence $ (\ga_j)_{j\in\N} \subset \wh G$ going to infinity $ \lim \limits_{j\to\iy}\noop{\va(\ga_j)}=0 $.
\item  For every $i \in \{1,...,r\}$ and any properly converging sequence  $\ol \ga=(\ga_j)_{j \in \N}\subset \GA_i$ whose limit set is contained in $S_{i-1} $  (taking a subsequence if necessary) and for the mappings $\tilde{\nu}_j=\tilde{\nu}_{\ol\ga,j}: CB(S_{i-1})\to \B(\H_i)$, one has
\begin{eqnarray*}
&\lim \limits_{j\to\iy }
\big\|
\va(\ga_j)-\tilde{\nu}_{j}\big(\va\res{S_{i-1}}\big)\big\|_{op}=0. 
\end{eqnarray*}
 \end{enumerate}
\end{theorem}

At the end of this article, an equivalent, but much simpler description of $C^*(G)$ will be given in Theorem \ref{thm genauer - sl2r}. 

\section{The unitary dual of $SL(2,\R)$}\label{unitary dual}
\subsection{Introduction of the operator $K_u$} \label{Ku}
Now, an operator $K_u$ which is needed in order to describe the unitary dual of $G$ will be introduced using the Knapp-Stein operator: \\
~\\
Define
$$C^{\infty}(K)_+:=\big\{f \in C^{\iy}(K)|~f(k)=f(-k)~\forall k \in K \big\}~~~\text{and}$$
$$C^{\infty}(K)_-:=\big\{f \in C^{\iy}(K)|~f(k)=-f(-k)~\forall k \in K \big\}$$
and let
$$J_u:C^{\infty}(K)_+ \to C^{\infty}(K)_+~~~\text{for}~u \in \C~\text{with}~\Re{u}>0$$
be  the Knapp-Stein intertwining operator, which sends the representation $\P^{+,u}$ to the representation $\P^{+,-u}$. Furthermore, let
$w:=k_{\frac{\pi}{2}}=\begin{pmatrix}
0 & -1  \\
1 & 0
\end{pmatrix}$ and extend $f \in L^2(K)$ to $G$ by using the Iwasawa decomposition $G \ni g=\kappa(g)e^{H(g)}\mu$ for $\kappa(g) \in K$, $H(g) \in \a$ and $\mu \in N$. Then, defining $\tilde{f}_u\big(\kappa(g) e^{H(g)}\mu \big):=e^{-(\nu_u+\rho)H(g)}f \big(\kappa(g)\big)$, this operator can be written as
\begin{eqnarray} \label{Integralformel}
J_uf(k)=\int \limits_N \tilde{f}_u(k \mu w) d \mu~~~\forall f \in C^{\infty}(K)_+~\forall k \in K.
\end{eqnarray}
This integral converges for $\Re{u}>0$ (see [\ref{knapp}], Chapter VII or [\ref{walII}], Chapter 10.1). \\
The mapping $f \mapsto J_u f$ is continuous and the family of operators $\{J_u|~u \in \C\}$ is holomorphic in $u$ for $\Re{u}>0$ with respect to appropriate topologies (see [\ref{knapp}], Chapter VII.7 or [\ref{walII}], Chapter 10.1). \\
For $u \in \R_{>0}$ the operator $J_u$ is self-adjoint with respect to the usual $L^2(K)$-scalar product. \\
Moreover, one can extend the function $u \mapsto J_u$ meromorphically to $\C$ (see [\ref{walII}], Chapter 10.1). Then, for every $u \in \C$ for which the operator $J_u$ is regular, it is an intertwining operator from $\P^{+,u}$ to $\P^{+,-u}$.

\begin{remark}\label{commute}
~\\
The operator $J_u$ 
commutes with the projections $p_n$ for all $n \in \N$ and for every $u \in \C$ for which $J_u$ is regular.
\end{remark}

This can be seen as follows: Since $J_u$ is a $G$-intertwining operator, it is a $K$-intertwining operator as well and therefore, it leaves $\H_{\P^{+,u}}(n)$, the $n$-th isotypic component of the representation $\P^{+,u}$ which was defined above, invariant. Hence, one can easily conclude that $p_n \circ J_u = J_u \circ p_n$ and the assertion follows. \\
%
~\\
One can now deduce that the operators $J_u$ have the property
\begin{eqnarray}\label{cn}
{J_u} \res{\H_{\P^{+,u}}(n)}=c_n(u) \cdot \id \res{\H_{\P^{+,u}}(n)}~~~\text{for all even}~ n \in \Z,
\end{eqnarray}
as an equality of meromorphic functions, where $c_n: \C \to \C$ is a meromorphic function 
for every even $n \in \N$. This follows from the above remark together with the fact that $\H_{\P^{+,u}}(n)$ is one-dimensional. \\
~\\
Using standard integral formulas (see [\ref{knapp}], Chapter 5.6), one gets by (\ref{Integralformel}) for $u=1$
$$J_1(f)~=~c \int \limits_K f(k) dk~~~\forall f \in C^{\iy}(K)_+$$
for a constant $c >0$. 
Therefore and since $\H_{\P^{+,u}}(n)= \C \cdot e^{-in \cdot}$ for every even $n \in \N$ which was stated above, one gets
\begin{eqnarray}\label{c_n(1)}
c_n(1)&&
\begin{cases}
\not=0~~~\text{for}~ n=0 \\
=0~~~\text{for}~ n \in \Z \setminus \{0\}~ \text{even.}
\end{cases}
\end{eqnarray}
~\\
For the convenience of the reader, an explicit formula for the quotients of the functions $c_n$ will now be given. However, this formula will not be used in this article. Instead, further necessary properties of the intertwining operator $J_u$ and the functions $c_n$ will be concluded from the irreducibility of the representations $\P^{\pm,u}$.

\begin{remark}\label{quotients}
~\\
The quotients of the functios $c_n$ can be given by
$$\frac{c_n(u)}{c_0(u)}~=~\frac{(u-1)(u-3)\cdots \big( u-(|n|-1)\big)}{(u+1)(u+3)\cdots \big( u+(|n|-1)\big)} \cdot (-1)^{\frac{n}{2}}~~~\text{for}~u \in \C ~\text{and all even}~ n \in \Z.$$
\end{remark}

This formula can be deduced from a formula for $c_n(u)$ in terms of Gamma functions which can be found in [\ref{cohn}], and the Gamma function recurrence formula. Here, one has to remark that the definition of $c_n(u)$ in [\ref{cohn}] differs by a sign from its definition in this article. \\
~\\
Next, it can be shown that
\begin{eqnarray}\label{cn3}
c_0(u) \not=0~~~\text{for}~u \in (0,1):
\end{eqnarray}
For this, assume that $c_0(u)=0$ for an element $u \in (0,1)$. Then, by (\ref{cn}), the operator $J_u$ has a non-zero kernel. Furthermore, the representation $\P^{+,u}$ is irreducible on $C^{\iy}(K)_+$ (see Remark \ref{irred}). As $\ker{J_u}$ is a closed invariant subspace with respect to the representation $\P^{+,u}$, the kernel of $J_u$ has to be the whole space $C^{\iy}(K)_+$ and hence the operator $J_u$ is identically zero. But there is no $u \in \C$ with $\Re{u}>0$ such that the operator $J_u$ is identically zero (see [\ref{walII}], Chapter 10.1). Therefore, $c_0(u) \not=0$ for every $u \in (0,1)$. \\
~\\
Thus, one can define
$$\tilde{J}_u:=\frac{1}{c_0(u)}~J_u$$
for $u \in \C$ as a meromorphic function.

\begin{lemma}\label{id}
~\\
$\tilde{J}_u$ is regular at $u=0$ and $\tilde{J}_0=\id$.
\end{lemma}

Proof: \\ 
First, one can observe that on $\H_{\P^{+,u}}(0)$ the operator $\tilde{J}_{u}$ is always equal to the identity and that in particular 
${\tilde{J}_{0}} {\res{\H_{\P^{+,0}}(0)}}$ also equals $\id \res{\H_{\P^{+,0}}(0)}$. \\
Now, it has to be shown that $\tilde{J}_u$ is regular at $u=0$. \\
Since the mapping $u \mapsto \tilde{J}_{u}$ is an operator-valued meromorphic function for $u \in \C$, one can represent it locally as a Laurent series in $0$ with finite principal part and operators as coefficients in the following way:
$$\tilde{J}_{u}=\sum \limits_{k=k_0}^{\iy} L_k u^k~~~\text{on}~C^{\iy}(K)_+$$
for operators $L_k$ going from $C^{\iy}(K)_+$ to $C^{\iy}(K)_+$ for $k\geq k_0$ and where $k_0 \in \Z$ is the smallest number such that $L_{k_0} \not=0$. \\
Moreover, this gives 
$$L_{k_0}=\lim \limits_{u \to 0} u^{-k_0} \tilde{J}_{u}.$$
Now, for the above desired regularity of $\tilde{J}_{u}$ at $u=0$, it has to be shown that $k_0 \geq 0$. \\
So, assume that $k_0 <0$. \\
As every $\tilde{J}_{u}$ is an intertwining operator from $\P^{+,u}$ to $\P^{+,-u}$, $L_{k_0}$ is an intertwining operator from $\P^{+,0}$ to itself, i.e. it commutes with $\P^{+,0}$. Furthermore, $L_{k_0}$ vanishes on $\H_{\P^{+,0}}(0)$, the space of all constant functions on $K$, because $\tilde{J}_{u}$ equals the identity on $\H_{\P^{+,u}}(0)$ and thus does not have a pole there. \\
Moreover, since $L_{k_0}$ commutes with $\P^{+,0}$, it va\-ni\-shes on $\P^{+,0}(g)\H_{\P^{+,0}}(0)$ for eve\-ry $g \in G$. In addition, the representation $\P^{+,0}$ is irreducible on the space $C^{\iy}(K)_{+}$ (see Remark \ref{irred}).  
Furthermore, for eve\-ry $ 0 \ne \xi \in \H_{\P^{+,0}}(0)$, the subspace $\text{span}\big\{ \P^{+,0}(g) \xi |~ g \in G \big\}$ is $G$-invariant and hence, by the irreducibility, $C^{\iy}(K)_{+} \subset \ol{\text{span} \big\{ \P^{+,0}(g) \xi |~ g \in G \big\}}$. Thus, $L_{k_0}$ va\-ni\-shes on the whole space $C^{\iy}(K)_{+}$, which is a contradiction to the choice of $k_0$. \\
Hence, one gets $k_0 \geq 0$, which means that the mapping $u \mapsto \tilde{J}_{u}$ does not have any poles in $u=0$. Therefore, $\tilde{J}_{u}$ is regular at $u=0$ on $C^{\iy}(K)_{+}$. \\
Moreover, as above, as a limit of intertwining operators, $\tilde{J}_{0}=\lim \limits_{u \to 0} \tilde{J}_{u}$ is an intertwining operator, which intertwines the irreducible unitary representation $\P^{+,0}$ with itself (see Remark \ref{irred} for the irreducibility and the unitarity of $\P^{+,0}$). Hence, by Schur's Lemma 
$\tilde{J}_{0}$ is a scalar multiple of the identity. Since it equals the identity on $\H_{\P^{+,0}}(0)$, one gets $\tilde{J}_{0}=\id$. \\
\qed
~\\
~\\
From (\ref{cn3}) and Lemma \ref{id}, one can conclude that
\begin{eqnarray}\label{cn4}
\frac{c_n(u)}{c_0(u)}>0~~~\text{for}~u \in (0,1)~\text{and for all even}~ n \in \Z:
\end{eqnarray}
Using the same argumentation for the operator $\tilde{J}_u$ as in the proof of (\ref{cn3}), one also gets that $\frac{c_n(u)}{c_0(u)} \not=0$ for every $u \in (0,1)$. \\
Furthermore, from Lemma \ref{id} one can deduce that $\frac{c_n(0)}{c_0(0)}=1$, and with the continuity of $\frac{c_n}{c_0}$ on $(0,1)$, one gets $\frac{c_n(u)}{c_0(u)}>0$ for all $u \in (0,1)$ and all even $n \in \Z$ and (\ref{cn4}) is shown.\\
~\\
Now, define a scalar product on 
$C^{\iy}(K)_+$ as follows:
$$\langle f_1,f_2 \rangle_u:=\big\langle \tilde{J}_u f_1,f_2 \big\rangle_{L^2(K)}.$$

\begin{lemma}\label{SPu}
~\\
$\langle \cdot,\cdot \rangle_u$ is an invariant positive definite scalar product for $u \in (0,1)$.
\end{lemma}
 
Proof: \\
$\langle \cdot, \cdot \rangle_u$ is hermitian: \\
By the definition of $\tilde{J}_u$ and as $J_u$ is self-adjoint with respect to the usual $L^2(K)$-scalar product for $u \in (0,1)$, this is straightforward. \\
$\langle \cdot, \cdot \rangle_u$ is invariant: \\
For every $g \in G$ the operator $\big(\P^{+,u}(g) \big)^{-1}$ is the adjoint operator of $\P^{+,-u}(g)$ with respect to the usual $L^2(K)$-scalar product. Now, let $f_1,f_2 \in C^{\iy}(K)_+$. Then, as $\tilde{J}_{u}$ intertwines $\P^{+,u}$ and $\P^{+,-u}$, one gets for every $g \in G$
\begin{eqnarray*}
\big\langle \P^{+,u}(g)f_1, \P^{+,u}(g)f_2 \big\rangle_{u} &=& \big\langle \tilde{J}_{u}\circ \P^{+,u}(g)f_1, \P^{+,u}(g)f_2 \big\rangle_{L^2(K)} \\
&=& \big\langle \P^{+,-u}(g)\circ \tilde{J}_{u}f_1, \P^{+,u}(g)f_2 \big\rangle_{L^2(K)} \\
&=& \big\langle \tilde{J}_{u}f_1, f_2 \big\rangle_{L^2(K)}~=~\langle f_1, f_2 \rangle_{u}.
\end{eqnarray*}
$\langle \cdot, \cdot \rangle_u$ is positive definite: \\
As $\frac{c_n(u)}{c_0(u)}>0$ by (\ref{cn4}), one gets for every $n \in \Z$ and $f \in \H_{\P^{+,u}}(n)$ 
by (\ref{cn}) above,
\begin{eqnarray*}
\langle f,f \rangle_u~=~ \big\langle \tilde{J}_u f,f \big\rangle_{L^2(K)}
~=~ \frac{c_n(u)}{c_0(u)}~ \Big\langle \id \res{ \H_{\P^{+,u}}(n)} f,f \Big\rangle_{L^2(K)}
~=~\frac{c_n(u)}{c_0(u)}~ \langle f,f \rangle_{L^2(K)} \geq 0 ~~~\text{and}
\end{eqnarray*}
\begin{eqnarray*}
\langle f,f \rangle_u=0 ~\Longleftrightarrow ~ \langle f,f \rangle_{L^2(K)}=0 ~\Longleftrightarrow~ f =0. 
\end{eqnarray*}
Now, since the direct sum $\bigoplus \limits_{n \in \Z} \H_{\P^{+,u}}(n)$
is orthogonal with respect to $\langle \cdot, \cdot \rangle_u$ and dense in $\H_{\P^{+,u}}=L^2(K)_+$, these observations hold for all $f \in C^{\iy}(K)_+$.\\
\qed 
~\\
~\\
~\\
The completion of $C^{\iy}(K)_+$ with respect to this scalar product $\langle \cdot, \cdot \rangle_u$ gives a Hilbert space $\H_u$. Considering the restriction of the representation $\P^{+,u}$ to $C^{\iy}(K)_+$ and then continuously extending it to the space $\H_u$, one gets a unitary representation which will be denoted by $\P^{+,u}$ as well. $G$ acts on $\H_u$ by this unitary representation $\P^{+,u}$. \\
~\\
Furthermore, let $d_n(u):=\sqrt{\frac{c_n(u)}{c_0(u)}}>0$ for $u \in (0,1)$. Next, a unitary bijection
$$K_u: \H_u \to L^2(K)_+~~~\forall u \in (0,1)$$
shall be defined. On the $n$-th isotypic component in $\H_u$, define $K_u$ by
$${K_u} \res{\H_{\P^{+,u}}(n)}:=d_n(u) \cdot \id \res{\H_{\P^{+,u}}(n)}~~~\text{for all even}~n \in \Z.$$
Then one can extend this definition to finite sums of $K$-types. 
This operator also is self-adjoint with respect to the usual $L^2(K)$-scalar product and for finite sums of $K$-types $f_1$ and $f_2$
$$\big\langle K_u f_1,K_u f_2 \big\rangle_{L^2(K)}=\big\langle K_u^2 f_1,f_2 \big\rangle_{L^2(K)}=\big\langle \tilde{J}_uf_1,f_2 \big\rangle_{L^2(K)}=\langle f_1,f_2 \rangle_u.$$
From this, it follows directly that it is unitary (if one regards the space $\H_u$ equipped with $\langle \cdot, \cdot \rangle_u$ and $L^2(K)_+$ with $\langle \cdot, \cdot \rangle_{L^2(K)}$) and hence, because of the density of $\bigoplus \limits_{n \in \Z} \H_{\P^{+,u}}(n)$ in $\H_u$, one can extend $K_u$ continuously on the whole space $\H_u$. \\
Moreover, $K_u$ is continuous in $u$ and one also has the property $\lim \limits_{u \to 0} K_u=\id$. \\
By its definition, the operator $K_u$ commutes with the projections $p_n$ for all $n \in \N$ as well. \\
~\\
Now, by [\ref{knapp}], Chapter 14.4, one gets the identity
\begin{eqnarray}\label{gam}
J_{-u} \circ J_u=\gamma(u) \cdot \id,
\end{eqnarray}
where $u \mapsto \gamma(u)$ is a meromorphic function. One can also obtain this equation by observing that $J_{-u} \circ J_u$ is an intertwining operator of the representation $\P^{+,u}$, which is irreducible for almost all $u \in \C$ (see Remark \ref{irred}), with itself, and by using a version of Schur's Lemma for $(\g,K)$-modules. \\
By restricting the above equation to the $n$-th isotypic component, one gets thus the relation
\begin{eqnarray}\label{cn-relation}
c_n(u)c_n(-u)=\gamma(u)~~~\text{for all even}~ n \in \Z
\end{eqnarray}
as meromorphic functions. \\
~\\
Next, another scalar product on the space $C^{\iy}(K)_{++}:=\big\{ f \in C^{\iy}(K)_{+}|~p_n(f)=0~\forall n \leq 0\big\}$ is needed: \\
For this, define
$$\tilde{J}_{(u)}:=\frac{1}{c_2(u)}~J_{u}$$ 
for $u \in \C$ as a meromorphic family of operators. Then, on $\H_{\P^{+,u}}(2)$ the operator $\tilde{J}_{(u)}$ is equal to the identity for every $u \in \C$ and in particular ${\tilde{J}_{(1)}} {\res{\H_{\P^{+,1}}(2)}}$ equals $\id \res{\H_{\P^{+,1}}(2)}$. \\
Moreover, define the space $C^{\iy}(K)_{+-}:=\big\{ f \in C^{\iy}(K)_{+}|~p_n(f)=0~\forall n \geq 0\big\}$. \\
Then, the representation $\P^{+,1}$ is irreducible on the spaces $C^{\iy}(K)_{++}$ and $C^{\iy}(K)_{+-}$ (see [\ref{wal}], Chapter 5.6).

\begin{lemma}\label{existence}
~\\
\begin{enumerate}[(a)]
\item $ \tilde{J}_{(-u)} \circ \tilde{J}_{(u)}= \tilde{J}_{(u)} \circ \tilde{J}_{(-u)}= \id $ as a meromorphic family.
\item $\tilde{J}_{(u)}$ is regular at $u=-1$ and 
$$\text{ker} \big( \tilde{J}_{(-1)}\big) \cap C^{\iy}(K)_{++}~=~\text{ker} \big( \tilde{J}_{(-1)}\big) \cap C^{\iy}(K)_{+-}~=~\{0\}.$$
Furthermore, $\tilde{J}_{(-1)}$ is an intertwining operator of $\P^{+,-1}$ with $\P^{+,1}$.
\item ${\tilde{J}_{(u)}} {\res{C^{\iy}(K)_{++} \oplus C^{\iy}(K)_{+-}}}$ is regular at $u=1$.
\end{enumerate}
\end{lemma}

Proof: \\
One has for all $n \in \Z$ and all $u \in \C$
$${\tilde{J}_{(-u)} \circ \tilde{J}_{(u)}} \res{\H_{\P^{+,u}}(n)} ~\overset{(\ref{gam})}{=}~\frac{\gamma(u)}{c_2(u)c_2(-u)} \cdot {\id} \res{\H_{\P^{+,u}}(n)} ~\overset{(\ref{cn-relation})}{=}~{\id} \res{\H_{\P^{+,u}}(n)}$$
and thus
$$\tilde{J}_{(-u)} \circ \tilde{J}_{(u)}~=~\tilde{J}_{(u)} \circ \tilde{J}_{(-u)}~=~ \id$$
for $u \in \C$ as a meromorphic family and (a) follows. \\
~\\
Because the mapping $u \mapsto \tilde{J}_{(u)}$ is an operator-valued meromorphic function for $u \in \C$, one can represent it locally as a Laurent series in $-1$ with finite principal part and operators as coefficients in the following way:
$$\tilde{J}_{(u)}=\sum \limits_{k=k_0}^{\iy} L_k (u+1)^k~~~\text{on}~C^{\iy}(K)_+$$
for operators $L_k$ going from $C^{\iy}(K)_+$ to $C^{\iy}(K)_+$ for $k\geq k_0$ and where $k_0 \in \Z$ is the smallest number such that $L_{k_0} \not=0$. \\
Then, one gets
\begin{eqnarray}\label{Lk0-Limes}
L_{k_0}~=~\lim \limits_{u \to 1} (-u+1)^{-k_0} \tilde{J}_{(-u)}~=~\lim \limits_{u \to 1} (-u+1)^{-k_0} \frac{J_{-u}}{c_2(-u)}
\end{eqnarray}
and thus, as a limit of intertwining operators, $L_{k_0}$ is an intertwining operator of $\P^{+,-1}$ with $\P^{+,1}$.
Now, it shall first be shown that
\begin{eqnarray}\label{ZwBeh}
\ker{L_{k_0}} \cap C^{\iy}(K)_{++}~=~\ker{L_{k_0}} \cap C^{\iy}(K)_{+-}~=~\{0\}. 
\end{eqnarray}
Since $\H_{\P^{+,1}}(n)=\H_{\P^{+,-1}}(n)$, 
$${L_{k_0}} \res{\H_{\P^{+,1}}(n)}=a_n \cdot {\id} \res{\H_{\P^{+,1}}(n)}~~~\text{for}~~~a_n:= \lim \limits_{u \to 1} \frac{c_n(-u)}{c_2(-u)}(-u+1)^{-k_0}.$$
As $C^{\iy}(K)_{++}= \bigoplus \limits_{\substack{n \in \Z_{>0} \\ n~even}} \H_{\P^{+,1}}(n)$ and $C^{\iy}(K)_{+-}= \bigoplus \limits_{\substack{n \in \Z_{<0} \\ n~even}} \H_{\P^{+,1}}(n)$, for (\ref{ZwBeh}) it has to be shown that $a_n \not=0$ for all even $n \in \Z \setminus \{0\}$. \\
By (\ref{cn-relation}), one has 
$$\frac{c_n(u)c_n(-u)}{c_2(-u)}~=~\frac{\gamma(u)}{c_2(-u)},$$
i.e. the left hand side does not depend on $n$. So the order of pole for the limit for $u \to 1$ has to be the same for every $n \in \Z$. But as $c_n(1)=0$ for $n \not=0$ and $c_0(1) \not= 0$ by (\ref{c_n(1)}), $c_0(-1)$ has to be $0$ in order to get the same order of pole for $n=0$ and $n \not=0$. It follows that $a_0=0$. As $k_0$ was chosen in such a way that $L_{k_0} \not=0$, there has to exist $0 \not=n \in \Z$ with $a_n \not=0$. \\
Furthermore, from (\ref{Integralformel}), one can conclude that $\overline{J_u f}=J_{\overline{u}} \overline{f}$ for all $u \in \C$ with $\Re{u}>0$. By extending meromorphically, this holds for all $u \in \C$ and, by regarding the $n$-th isotypic component, one can deduce that $\overline{c_n(u)}=c_{-n}(\overline{u})$.  
Hence, one gets for every $n \in \Z$ and every $u \in \R$ that $\overline{c_n(u)}=c_{-n}(u)$. This means that there also exists an integer $n_1 >0$ such that $a_{n_1} \not=0 \not= a_{-n_1}$. Therefore,
$${L_{k_0}} \res{\H_{\P^{+,1}}(n_1)}~=~a_{n_1} \cdot {\id} \res{\H_{\P^{+,1}}(n_1)} ~\not=0 \not=~ a_{-n_1} \cdot {\id} \res{\H_{\P^{+,1}}(-n_1)}~= ~{L_{k_0}}\res{\H_{\P^{+,1}}(-n_1)}$$
and thus
$$\im{ \Big({L_{k_0}}\res{\H_{\P^{+,1}}(n_1)} \Big)} \subset \H_{\P^{+,1}}(n_1) \subset C^{\iy}(K)_{++}~~~\text{and}$$
$$\im{ \Big({L_{k_0}}\res{\H_{\P^{+,1}}(-n_1)} \Big)} \subset \H_{\P^{+,1}}(-n_1) \subset C^{\iy}(K)_{+-}.$$
So, in particular, $\Im{(L_{k_0})} \cap C^{\iy}(K)_{++} \not=\{0\}\not= \Im{(L_{k_0})} \cap C^{\iy}(K)_{+-}$. But as $\P^{+,1}$ is irreducible on $C^{\iy}(K)_{++}$ and on $C^{\iy}(K)_{+-}$, one gets
$$\ol{\Im{(L_{k_0})} \cap C^{\iy}(K)_{++}}~=~C^{\iy}(K)_{++}~~~\text{and}~~~\ol{\Im{(L_{k_0})} \cap C^{\iy}(K)_{+-}}~=~C^{\iy}(K)_{+-},$$
where the completions are regarded in the spaces $C^{\iy}(K)_{++}$ and respectively $C^{\iy}(K)_{+-}$. Hence, $\H_{\P^{+,1}}(n) \subset \Im{(L_{k_0})}$ for all even $n \not= 0$ and therefore, $a_n \not=0$ for every even $n \not= 0$ and (\ref{ZwBeh}) follows. \\
One can now conclude that $k_0=0$: \\
Since
$${L_{k_0}}\res{\H_{\P^{+,1}}(2)}~= ~\lim \limits_{u \to 1} (-u+1)^{-k_0} \cdot {\id}\res{\H_{\P^{+,1}}(2)}$$
does not have any poles, $k_0 \not> 0$, and as it does not have any zeros, $k_0 \not< 0$. \\
Hence, 
$$L_{k_0}= \lim \limits_{u \to 1} \tilde{J}_{(-u)}=\tilde{J}_{(-1)}$$
and thus, (b) follows by (\ref{ZwBeh}) and as $L_{k_0}$ intertwines $\P^{+,-1}$ with $\P^{+,1}$. \\
Now, by (b), $\tilde{J}_{(u)}$ is regular at $u=-1$. As by (a) one has $\tilde{J}_{(-1)} \circ \tilde{J}_{(1)}= \tilde{J}_{(1)} \circ \tilde{J}_{(-1)}= \id $ and since by (b) the operator $\tilde{J}_{(-1)}$ is nowhere equal to $0$ on $C^{\iy}(K)_{++} \oplus C^{\iy}(K)_{+-}$, the operator $\tilde{J}_{(u)}$ has to be regular on $C^{\iy}(K)_{++} \oplus C^{\iy}(K)_{+-}$ at $u=1$ as well and (c) is shown. \\
\qed 
~\\
~\\
Moreover, $\tilde{J}_{(1)}$ is not identically $0$ on the space $C^{\iy}(K)_{++}$, as it equals the identity on $\H_{\P^{+,1}}(2)$. Thus, regard the operator $\tilde{J}_{(1)}$ on the space $C^{\iy}(K)_{++}$. This operator is injective, since for eve\-ry $f \in C^{\iy}(K)_{++}$ with $\tilde{J}_{(1)}(f)=0$, one gets by Lemma \ref{existence}(a) that $0=\tilde{J}_{(-1)} \circ \tilde{J}_{(1)}(f)=\id(f)=f$. \\
~\\
Now, define for all functions $f_1, f_2 \in C^{\iy}(K)_{++}$
$$\langle f_1,f_2 \rangle_{(1)}:= \big\langle \tilde{J}_{(1)} f_1,f_2 \big\rangle_{L^2(K)}.$$
Furthermore, choose $\tilde{f_1}$ in such a way that $\tilde{J}_{(-1)}\big(\tilde{f_1} \big)=f_1$. 
Then
\begin{eqnarray}\label{sp-gleichheit}
\langle f_1,f_2 \rangle_{(1)}=\big\langle \tilde{f_1},f_2 \big\rangle_{L^2(K)},
\end{eqnarray}
since the scalar product does not depend on the choice of $\tilde{f_1}$: \\
It is possible to add an element $\tilde{f}_{ker} \in \text{ker} \big(\tilde{J}_{(-1)}\big)$ to the function $\tilde{f_1}$. But by Lemma \ref{existence}(b) above, one has $\text{ker} \big( \tilde{J}_{(-1)}\big) \cap C^{\iy}(K)_{++}=\text{ker} \big( \tilde{J}_{(-1)}\big) \cap C^{\iy}(K)_{+-}=\{0\}$ and therefore, \mbox{$\tilde{f}_{ker} \in \big\{ f \in C^{\iy}(K)_{+}|~p_n(f)=0~\forall n \not= 0\big\}$}. Thus, the function $\tilde{f}_{ker}$ is orthogonal to \mbox{$f_2 \in C^{\iy}(K)_{++}=\big\{ f \in C^{\iy}(K)_{+}|~p_n(f)=0~\forall n \leq 0\big\}$} and the scalar product stays the same. 
 
\begin{lemma}\label{(1)-SP}
~\\
$\langle \cdot,\cdot \rangle_{(1)}$ is an invariant positive definite scalar product.
\end{lemma}
 
Proof: \\
$\langle \cdot, \cdot \rangle_{(1)}$ is hermitian: \\
As above for $\langle \cdot, \cdot \rangle_u$, this is straightforward. \\
$\langle \cdot, \cdot \rangle_{(1)}$ is invariant: \\
Let $f_1,f_2 \in C^{\iy}(K)_{++}$ and choose $\tilde{f_1}$ such that $\tilde{J}_{(-1)}\big(\tilde{f_1} \big)=f_1$. \\
Then, as $\tilde{J}_{(-1)}$ intertwines $\P^{+,-1}$ and $\P^{+,1}$ by Lemma \ref{existence}(b), for every $g \in G$
$$\tilde{J}_{(-1)} \Big( \P^{+,-1}(g)  \tilde{f_1}  \Big)~=~\P^{+,1}(g) \circ \tilde{J}_{(-1)} \big(\tilde{f_1} \big)~=~\P^{+,1}(g) f_1.$$
Hence, one can choose 
$$\widetilde{\P^{+,1}(g) f_1}:=\P^{+,-1}(g)  \tilde{f_1}.$$
Since for all $g \in G$ the operator $\big(\P^{+,1}(g) \big)^{-1}$ is the adjoint operator of $\P^{+,-1}(g)$ with respect to the usual $L^2(K)$-scalar product, one gets for every $g \in G$
\begin{eqnarray*}
\big\langle \P^{+,1}(g)f_1, \P^{+,1}(g)f_2 \big\rangle_{(1)} &=& \Big\langle \widetilde{\P^{+,1}(g)f_1}, \P^{+,1}(g)f_2 \Big\rangle_{L^2(K)}\\ 
&=& \Big\langle \P^{+,-1}(g)  \tilde{f_1} , \P^{+,1}(g)f_2 \Big\rangle_{L^2(K)} \\
&=& \big\langle \tilde{f_1} , f_2 \big\rangle_{L^2(K)} \\
&=& \langle f_1 , f_2 \rangle_{(1)}. 
\end{eqnarray*}
$\langle \cdot, \cdot \rangle_{(1)}$ is positive definite: \\
$\frac{c_n(u)}{c_2(u)}>0$ for every $u \in (0,1)$, since
$$\frac{c_n(u)}{c_2(u)}=\frac{c_n(u)}{c_0(u)} \cdot \frac{c_0(u)}{c_2(u)}>0$$
by (\ref{cn4}) of the beginning of this subsection. Hence, its limit $\lim \limits_{u \to 1} \frac{c_n(u)}{c_2(u)}$ is larger or equal to $0$ as well. Therefore, similar as above for $\langle \cdot, \cdot \rangle_u$, for every $n \in \N^*:= \N \setminus \{0\}$ and for every $f_1 \in \H_{\P^{+,1}}(n)$, one gets
\begin{eqnarray*}
\langle f_1,f_1 \rangle_{(1)}&=& \big\langle \tilde{J}_{(1)} f_1,f_1 \big\rangle_{L^2(K)} 
~=~ \lim \limits_{u \to 1} \frac{c_n(u)}{c_2(u)}~ \langle f_1,f_1 \rangle_{L^2(K)} \geq 0
\end{eqnarray*}
from the argument above and since $\H_{\P^{+,1}}(n)=\H_{\P^{+,u}}(n)$ for every $u \in (0,1)$. Moreover, because of the injectivity of $\tilde{J}_{(1)}$
\begin{eqnarray*}
 f_1=0 ~\Longleftrightarrow ~\tilde{J}_{(1)}  f_1=0 ~\Longleftrightarrow~  \lim \limits_{u \to 1} \frac{c_n(u)}{c_2(u)}~ f_1=0,
\end{eqnarray*}
which means that $\lim \limits_{u \to 1} \frac{c_n(u)}{c_2(u)}>0$. Hence,
\begin{eqnarray*}
\langle f_1,f_1 \rangle_{(1)}=0 ~\Longleftrightarrow ~ \langle f_1,f_1 \rangle_{L^2(K)}=0 ~\Longleftrightarrow~ f_1 =0. 
\end{eqnarray*}
As the direct sum $\bigoplus \limits_{n \in \N^*} \H_{\P^{+,1}}(n)$ is orthogonal with respect to $\langle \cdot, \cdot \rangle_{(1)}$ and dense in $C^{\iy}(K)_{++}$, the positive definiteness is shown everywhere. \\
\qed
~\\
~\\
Now, the completion of the space $C^{\iy}(K)_{++}$ with respect to this scalar product gives a Hilbert space which will be called $\H_{(1)}$. \\
~\\
The same procedure can be accomplished for $C^{\iy}(K)_{+-}=\big\{ f \in C^{\iy}(K)_{+}|~p_n(f)=0~\forall n \geq 0\big\}$: \\
Here, for $u \in (0,1)$, define the operator $\tilde{J}_{[u]}$ as
$$\tilde{J}_{[u]}:=\frac{1}{c_{-2}(u)}~J_{u}.$$
As above, it can be shown that $\tilde{J}_{[u]}$ is regular at $u=1$ on $C^{\iy}(K)_{+-}$ and is not identically $0$ there. \\
Thus, define again for all $f_1, f_2 \in C^{\iy}(K)_{+-}$
$$\langle f_1,f_2 \rangle_{[1]}:= \big\langle \tilde{J}_{[1]} f_1,f_2 \big\rangle_{L^2(K)}.$$
This is another invariant positive definite scalar product. The completion of the space $C^{\iy}(K)_{+-}$ with respect to this scalar product gives a Hilbert space called $\H_{[1]}$.

\subsection{Description of the irreducible unitary representations}\label{rep}
Now, some convenient realizations for the unitary dual of $SL(2,\R)$ shall be provided. \\
~\\
The unitary dual $\wh{G}$ of $G=SL(2,\R)$ consists of the following representations:

\begin{enumerate}
\item The \textbf{principal series} representations:
\begin{enumerate}
\item $\P^{+,iv}$ for $v \in [0, \iy)$. 
\item  $\P^{-,iv}$ for $v \in (0, \iy)$. 
\end{enumerate}
See Section \ref{prel-sl2r} above for the definitions.
\item The \textbf{complementary series} representations ${\cal C}^{u}$ for $u \in (0,1)$: \\
The Hilbert space $\H_{{\cal C}^{u}}$ is defined by
$$\H_{{\cal C}^{u}}:=L^2(K)_+$$
and the action is given by 
$${\cal C}^{u}(g):=K_u \circ \P^{+,u} (g)\circ K_u^{-1}$$
for all $g \in G$, where here again $\P^{+,u} (g)$ is meant in the following way: One considers the restriction of $\P^{+,u} (g)$ to $C^{\iy}(K)_+$ and then continuously extends it to the space $\H_u$ (see the definition of $\H_u$ in Section \ref{Ku}).
\item The \textbf{discrete series} representations:
\begin{enumerate}
\item $\D_{m}^{+}$ for odd $m \in \N^*$: 
\begin{enumerate}
\item $\D_{1}^+$: \\
The Hilbert space $\H_{\D_1^+}$ is given by
$$\H_{\D_1^+}:=\H_{(1)}$$
defined in Section \ref{Ku}. 
The action is given by 
$$\D_1^+:=\P^{+,1}.$$
Here again, as well as in all the definitions in this subsection, the representation $\P^{+,u}$ for the different values $u \in \C$ is meant as described above: One restricts it to the respective subspace of $L^2(K)$ and then continuously extends it to the respective Hilbert space. 
\item $\D_{m}^{+}$ for odd $m \in \N_{\geq 3}$: \\
As a Hilbert space $\H_{\D_m^+}$ for $\D_{m}^+$ for odd $m \in \N_{\geq 3}$ one can take the completion of the space
$$\big\{f \in C^{\iy}(K)_+|~p_n(f)=0~\forall n \leq m-1\big\}$$
with respect to an appropriate scalar product, and as the action 
$$\D_m^+:=\P^{+,m}.$$
With this realization, the Hilbert spaces $\H_{\D_{m}^{+}}$ for odd $m \in \N_{\geq 3}$ depend on $m$. But as all of them are infinite-dimensional and separable, one can identify them if one conjugates the respective $G$-action. So, fix an infinite-dimensional separable Hilbert space $\H_{\D}$. Furthermore, the $G$-action is not needed for the determination of $C^*(G)$. 
\end{enumerate}
\item $\D_{m}^{-}$ for odd $m \in \N^*$: 
\begin{enumerate}
\item $\D_{1}^-$: \\
The Hilbert space $\H_{\D_1^-}$ is given by
$$\H_{\D_1^-}:=\H_{[1]}$$
defined in Section \ref{Ku} above and the action is given by 
$$\D_1^-:=\P^{+,1}.$$
\item $\D_{m}^-$ for odd $m \in \N_{\geq 3}$: \\
Similar as for $\D_m^+$, as a Hilbert space $\H_{\D_m^-}$ for $\D_{m}^-$ for odd $m \in \N_{\geq 3}$ one can take the completion of the space
$$\big\{f \in C^{\iy}(K)_+|~p_n(f)=0~\forall n \geq -m+1\big\}$$
with respect to an appropriate scalar product, and as the action 
$$\D_m^-:=\P^{+,m}.$$
Again, the Hilbert spaces depend on $m$. One identifies them and takes the common infinite-dimensional separable Hilbert space $\H_{\D}$ fixed in (a)(ii). Again, the $G$-action is not needed for the determination of $C^*(G)$.
\end{enumerate}
\item $\D_{m}^{+}$ for even $m \in \N^*$: \\
As a Hilbert space $\H_{\D_m^+}$ for $\D_{m}^+$ for even $m \in \N^*$ one can take the completion of the space
$$\big\{f \in C^{\iy}(K)_-|~p_n(f)=0~\forall n \leq m-1\big\}$$
with respect to an appropriate scalar product, and as the action 
$$\D_m^+:=\P^{-,m}.$$
Again, the Hilbert spaces are identified and one takes the common infinite-dimensional se\-parable Hilbert space $\H_{\D}$, as in (a)(ii).
\item $\D_{m}^-$ for even $m \in \N^*$: \\
As a Hilbert space $\H_{\D_m^-}$ for $\D_{m}^-$ for even $m \in \N^*$ one can take the completion of the space
$$\big\{f \in C^{\iy}(K)_-|~p_n(f)=0~\forall n \geq -m+1\big\}$$
with respect to an appropriate scalar product, and as the action 
$$\D_m^-:=\P^{-,m}.$$
Here again, the Hilbert spaces are identified and one takes the common Hilbert space $\H_{\D}$, as in (a)(ii).
\end{enumerate}
\item The \textbf{limits of the discrete series}:
\begin{enumerate}
\item $\D_+$: \\
The Hilbert space $\H_{\D_+}$ is defined by
$$\H_{\D_+}:=\big\{f \in L^2(K)_-|~p_n(f)=0~\forall n \leq 0 \big\}$$
and the action is given by 
$$\D_+:=\P^{-,0}.$$
\item $\D_-$: \\
The Hilbert space $\H_{\D_-}$ is defined by
$$\H_{\D_-}:=\big\{f \in L^2(K)_-|~p_n(f)=0~\forall n \geq 0 \big\}$$
and the action is given by 
$$\D_-:=\P^{-,0}.$$
\end{enumerate}
\item The \textbf{trivial} representation $\F_1$: \\
Its Hilbert space $\H_{\F_1}=\C$ will be identified with the space of constant functions
$$\big\{f \in L^2(K)_+|~p_n(f)=0~\forall n \not=0 \big\}.$$
Here, the action is given by
$$\F_1(g):=\id$$
for all $g \in G$. One also has
$$\F_1=\P^{+,-1},$$
as $\nu_{-1}+\rho=0$ and every $f \in \H_{\F_1}$ is a constant function. 
\end{enumerate}

A discussion of all irreducible representations of $SL(2, \R)$ without scalar products can be found in [\ref{wal}], Chapter 5.6. But as the scalar products in this article are shown to be unitary and invariant, they are the correct ones. \\
In [\ref{knapp}], Chapter 2.5, a different realization of the discrete series representations can be read that is not used in this article. An alternative description of all irreducible unitary representations can be found in [\ref{lang}], Chapter 6.6.

\begin{remark} \label{rem1}
$$\H_{\P^{-,0}} \cong \H_{\D_+} \oplus \H_{\D_-}.$$
\end{remark}

\begin{remark} \label{rem-pn}
~\\
By the definition of the operator $K_u$ by means of its value on the space $\H_{\P^{+,u}}$, one can easily verify that $p_n$ also projects $\H_{{\cal C}^u}$ to its $n$-th isotypic component $\H_{{\cal C}^u}(n)$ for every $n \in \Z$. \\
Furthermore, for every irreducible unitary representation $\pi$ of $G$ the operator $p_n$ leaves the above defined Hilbert space $\H_{\pi}$ invariant, since all of the Hilbert spaces are completions of the space $C^{\iy}(K)$ fulfilling $p_n$-cancelation properties for certain $n \in \Z$. \\
Hence, for every $n \in \Z$ and for every irreducible unitary representation $\pi$ of $G$, the operator $p_n$ is the projection of $\H_{\pi}$ to the $n$-th isotypic component $\H_{\pi}(n)$. 
\end{remark}

\subsection{The $SL(2, \R)$-representations applied to the Casimir operator} \label{Casi}

Regard the Casimir operator $\mathscr{C}$ in the universal enveloping algebra $\U(\g)$ with respect to the non-degenerate symmetric and $Ad$-invariant bilinear form on $\g$ defined by
$$\langle X, Y \rangle:= 2 \tr (XY)~~~\forall X,Y \in \g.$$

In order to be able to describe the topology on $\wh{G}$, the above listed representations applied to the Casimir operator will now be given:

\begin{lemma}\label{pcasi}
~\\
Applied to the Casimir operator, the representations $\P^{\pm,u}$ for $u \in \C$ give the following:
$$\P^{\pm,u}(\mathscr{C})=\frac{1}{4}\big(u^2-1 \big) \cdot \id.$$
\end{lemma}

The proof of this Lemma is standard and consists of an easy calculation using the fact that $\P^{\pm,u}(g) \circ \P^{\pm,u}(\mathscr{C})=\P^{\pm,u}(\mathscr{C}) \circ \P^{\pm,u}(g)$ for all $g \in G$. (Compare [\ref{wal}], Chapter 5.6.) \\
~\\
From Lemma \ref{pcasi} one can deduce for the above listed irreducible unitary representations of $SL(2, \R)$: 
\begin{enumerate}
\item $\P^{\pm,iv}(\mathscr{C})=\frac{1}{4} \big(-v^2-1 \big) \cdot \id~~~\forall v \in [0, \iy)$, \\
\item ${\cal C}^u(\mathscr{C})=K_u \circ \frac{1}{4}\big(u^2-1 \big)\cdot \id \circ K_u^{-1}=\frac{1}{4} \big(u^2-1 \big)\cdot \id~~~\forall u \in (0,1)$, \\
\item $\D_m^{\pm}(\mathscr{C})=\frac{1}{4} \big(m^2-1 \big)\cdot \id~~~\forall m \in \N^*$, \\
\item $\D_{\pm}(\mathscr{C})=\frac{1}{4}(0-1)\cdot \id=-\frac{1}{4}\cdot \id~~$ and \\
\item $\F_1(\mathscr{C})=\frac{1}{4} \big((-1)^2-1 \big)\cdot \id=0$.
\end{enumerate}

\subsection{The topology on $\wh{SL(2,\R)}$}\label{top}
With the help of the above computations, it is now possible to describe the topology on $\wh{G}$:

\begin{proposition}
~\\
The topology on $\wh{G}$ can be characterized in the following way: 
\begin{enumerate}
\item For all sequences $(v_j)_{j \in \N}$ and all $v$ in $[0, \iy)$
$$\P^{+,iv_j} \overset{j \to \iy}{\longrightarrow} \P^{+,iv}~ \Longleftrightarrow ~v_j \overset{j \to \iy}{\longrightarrow}v.$$ 
\item For all sequences $(v_j)_{j \in \N}$ and all $v$ in $(0, \iy)$
$$\P^{-,iv_j} \overset{j \to \iy}{\longrightarrow} \P^{-,iv} ~\Longleftrightarrow ~v_j \overset{j \to \iy}{\longrightarrow}v.$$ 
For all sequences $(v_j)_{j \in \N}$ in $[0, \iy)$
$$\P^{-,iv_j} \overset{j \to \iy}{\longrightarrow} \big\{\D_+,\D_-\big\} ~\Longleftrightarrow ~v_j \overset{j \to \iy}{\longrightarrow}0.$$
\item For all sequences $(u_j)_{j \in \N}$ and all $u$ in $(0,1)$
$${\cal C}^{u_j} \overset{j \to \iy}{\longrightarrow} {\cal C}^{u} ~\Longleftrightarrow~ u_j \overset{j \to \iy}{\longrightarrow}u.$$
For all sequences $(u_j)_{j \in \N}$ in $(0,1)$
$${\cal C}^{u_j} \overset{j \to \iy}{\longrightarrow} \P^{+,0} ~\Longleftrightarrow~ u_j \overset{j \to \iy}{\longrightarrow}0.$$
For all sequences $(u_j)_{j \in \N}$ in $(0,1)$
$${\cal C}^{u_j} \overset{j \to \iy}{\longrightarrow} \big\{\D_1^+,\D_1^-,\F_1\big\}~ \Longleftrightarrow ~ u_j \overset{j \to \iy}{\longrightarrow}1.$$
\end{enumerate}
All other sequences $(\pi_j)_{j \in \N}$ can only converge if they fulfill one of the following conditions: 
\begin{enumerate}[(a)]
\item They are constant for large $j \in \N$. In that case they have one single limit point, namely the value they take for large $j$. 
\item For large $j \in \N$ they are one of the above listed sequences. Then they converge in the above described way. 
\item For large $j \in \N$ the sequence consists only of members $\P^{+,iv_j}$ and ${\cal C}^{u_j}$. Then it converges to $\P^{+,0}$ if and only if $v_j \overset{j \to \iy}{\longrightarrow}0$ as well as $u_j \overset{j \to \iy}{\longrightarrow}0$ (compare 1. for $v=0$ and the second part of 3. above).
\item For large $j \in \N$ the sequence consists only of members ${\cal C}^{u_j}$ and $\D_1^+$. Then it converges to $\D_1^+$ if and only if $u_j \overset{j \to \iy}{\longrightarrow}1$ (compare the third part of 3. above).
\item For large $j \in \N$ the sequence consists only of members ${\cal C}^{u_j}$ and $\D_1^-$. Then it converges to $\D_1^-$ if and only if $u_j \overset{j \to \iy}{\longrightarrow}1$ (compare the third part of 3. above).
\item For large $j \in \N$ the sequence consists only of members ${\cal C}^{u_j}$ and $\F_1$. Then it converges to $\F_1$ if and only if $u_j \overset{j \to \iy}{\longrightarrow}1$ (compare the third part of 3. above).
\item For large $j \in \N$ the sequence consists only of members $\P^{-,iv_j}$ and $\D_+$. Then it converges to $\D_+$ if and only if $v_j \overset{j \to \iy}{\longrightarrow}0$ (compare the second part of 2. above).
\item For large $j \in \N$ the sequence consists only of members $\P^{-,iv_j}$ and $\D_-$. Then it converges to $\D_-$ if and only if $v_j \overset{j \to \iy}{\longrightarrow}0$ (compare the second part of 2. above).
\end{enumerate}
\end{proposition}

Since the topology of $\wh{G}$ for $G=SL(2, \R)$ is well-known, the proof of this proposition will be a bit sketchy, containing the ideas but leaving out most of the calculations: \\
If a sequence of representations $(\pi_j)_{j \in \N}$ converges to a representation $\pi$, the sequence $\big(\pi_j(\mathscr{C})\big)_{j \in \N}$ has to converge to $\pi(\mathscr{C})$ as well. By the observations of Section \ref{Casi}, the left hand side thus implies the right hand side in all the cases. \\
Now, for the other implication, it will first be shown, that for a sequence $(u_j)_{j \in \N}$ in $\C$, $u \in \C$ and for every compact set $\tilde{K} \subset G$
\begin{eqnarray} \label{hilfskonvergenz}
u_j \overset{j \to \iy}{\longrightarrow} u ~\Longrightarrow~ \sup \limits_{g \in \tilde{K}} \big\|\P^{\pm,u_j}(g)- \P^{\pm,u}(g) \big\|_{op} \overset{j \to \iy}{\longrightarrow} 0:
\end{eqnarray}
One has 
\begin{eqnarray*} 
&&\sup \limits_{g \in \tilde{K}} \big\|\P^{\pm,u_j}(g)- \P^{\pm,u}(g) \big\|_{op}^2 \\
&=&\sup \limits_{g \in \tilde{K}} ~\sup \limits_{\substack{f \in L^2(K)_{\pm}\\ \|f\|_2=1}}~ \int \limits_K \big| \P^{\pm,u_j}(g)(f)(k)- \P^{\pm,u}(g)(f)(k)\big|^2dk \\
&=&\sup \limits_{g \in \tilde{K}}~ \sup \limits_{\substack{f \in L^2(K)_{\pm}\\ \|f\|_2=1}}~
 \int \limits_K \Big|e^{-\nu_{u_j} H(g^{-1}k)}-e^{-\nu_{u} H(g^{-1}k)}\Big| ^2 e^{-2 \rho H(g^{-1}k)} \Big| f \Big(\kappa \big(g^{-1}k \big) \Big)\Big|^2 dk.
\end{eqnarray*}
Now, for $g \in \tilde{K}$ and $k \in K$, $g^{-1}k$ is also contained in a compact set and therefore, $H \big(g^{-1}k \big)$ is contained in a compact set as well. As
$$H \big(g^{-1}k \big)~=~
\begin{pmatrix}
h(g^{-1}k) & 0  \\
0 & -h(g^{-1}k)
\end{pmatrix}~~~\text{for}~h(g^{-1}k) \in \R,$$
there is thus a compact set $I \subset \R$ such that $h(g^{-1}k) \in I$ for all $g \in \tilde{K}$ and all $k \in K$. Hence
\begin{eqnarray*}
&&\sup \limits_{g \in \tilde{K}}~ \sup \limits_{\substack{f \in L^2(K)_{\pm}\\ \|f\|_2=1}}~
\int \limits_K \Big|e^{-\nu_{u_j} H(g^{-1}k)}-e^{-\nu_{u} H(g^{-1}k)}\Big| ^2 e^{-2 \rho H(g^{-1}k)} \Big| f \Big(\kappa \big(g^{-1}k \big) \Big)\Big|^2 dk \\
&\leq& \sup \limits_{x \in I} \Big|e^{-u_j x}-e^{-u x}\Big|^2
~\sup \limits_{g \in \tilde{K}} ~\sup \limits_{\substack{f \in L^2(K)_{\pm}\\ \|f\|_2=1}}~ \int \limits_K e^{-2 \rho H(g^{-1}k)}\Big|f \Big(\kappa \big(g^{-1}k \big) \Big) \Big|^2 dk \\
&\overset{Lemma~\ref{intabsch}}{=}& \sup \limits_{x \in I} \Big|e^{-u_j x}-e^{-u x}\Big|^2
~ \sup \limits_{\substack{f \in L^2(K)_{\pm}\\ \|f\|_2=1}} \| f\|_{L^2(K)}^2 \\
&=& \sup \limits_{x \in I} \Big|e^{-u_j x}-e^{-u x}\Big|^2 \overset{j \to \iy}{\longrightarrow}0. 
\end{eqnarray*}
Therefore, (\ref{hilfskonvergenz}) follows. \\
From this, one can easily deduce that for a sequence $(u_j)_{j \in \N}$ in $\C$ and $u \in \C$
\begin{eqnarray*}
u_j \overset{j \to \iy}{\longrightarrow} u ~\Longrightarrow ~ \P^{\pm,u_j} \overset{j \to \iy}{\longrightarrow} \P^{\pm,u}
\end{eqnarray*}
in the sense of convergence of matrix coefficients described in Theorem \ref{top-dual}, i.e. for some $f \in \H_{\P^{\pm,u}}$ there exists for every $j \in \N$ a function $f_j \in \H_{\P^{\pm,u_j}}$ such that
\begin{eqnarray}\label{conv}
\big\langle \P^{\pm,u_j}(\cdot) f_j,f_j \big\rangle_{\H_{\P^{\pm,u_j}}} \overset{j \to \iy}{\longrightarrow} \big\langle \P^{\pm,u}(\cdot) f,f \big\rangle_{\H_{\P^{\pm,u}}}
\end{eqnarray}
uniformly on compacta. \\
This yields directly the second implication of 1. and 2. (using that $\H_{\D_+}, \H_{\D_-} \subset \H_{\P^{-,iv}}$ for all $v \in (0, \iy)$). \\
Furthermore, one gets the equality
\begin{eqnarray}\label{intertwining}
\big\langle K_u \circ \P^{+,u}(g) \circ K_u ^{-1} f,f \big\rangle_{L^2(K)}~=~\big\langle \P^{+,u}(g) f,f \big\rangle_{L^2(K)}
\end{eqnarray}
for $\tilde{n} \in \Z$, $u \in (0,1)$, $g \in G$ and $f \in \H_{\P^{+,u}}(\tilde{n})$. \\
In order to show this, one uses the fact that the operator $K_u$ is self-adjoint with respect to the usual $L^2(K)$-scalar product. \\
Now, by (\ref{intertwining}), choosing $f$ in the space $\H_{\P^{\pm,u}}(n)$ for the matching value $n \in \Z$, one can express the matrix coefficients $\big\langle \pi(g) f,f \big\rangle_{\H_{\pi}}$ for every representation $\pi \in \wh{G}$ as $\big\langle \P^{+,u}(g) f,f \big\rangle_{L^2(K)}$. Then, one gets by (\ref{conv}) for $f_j=f$ the convergences needed for the second implications of 3., using Lemma \ref{existence} and the characterization of the $\H_{(1)}$-scalar product in (\ref{sp-gleichheit}). \\
At the end it still has to be shown that these are all possibilities of convergence. \\
For this, first, one can see that for a sequence of representations $(\pi_j)_{j \in \N}$ converging to a representation $\pi$, there has to be a common $K$-type for $(\pi_j)_{j~\text{large}}$ and $\pi$. With this fact and by considering the values the representations of $SL(2, \R)$ take on the Casimir operator (see Section \ref{Casi}), one can separate several sets of representations from each other such that only the possibilities of convergence listed in this proposition are possible. \\
\qed

\subsection{Definition of subsets $\Gamma_i$ of the unitary dual}\label{gammas}
Now, the unitary dual will be divided into different subsets which are thereafter proved to meet the requirements of the norm controlled dual limit conditions: \\
~\\
Define
\begin{eqnarray*}
\Gamma_0&:=&\big\{\F_1 \big\} \\
\Gamma_1&:=&\big\{\D_1^+ \big\} \\
\Gamma_2&:=&\big\{\D_1^- \big\} \\
\Gamma_3&:=&\big\{\D_+ \big\} \\
\Gamma_4&:=&\big\{\D_- \big\} \\
\Gamma_5&:=&\big\{\D_m^{\pm}|~m \in \N_{> 1} \big\} \\
\Gamma_6&:=&\big\{\P^{+,iv}|~v \in [0, \iy) \big\} \\
\Gamma_{7}&:=&\big\{\P^{-,iv}|~v \in (0, \iy) \big\} \\
\Gamma_{8}&:=&\big\{{\cal C}^{u}|~u \in (0,1) \big\}. 
\end{eqnarray*}

Obviously, all the sets $\Gamma_i$ for $i \in \{0,...,8\}$ are Hausdorff. Furthermore, the sets
$$S_i:=\bigcup \limits_{j \in \{0,...,i\}} \Gamma_j$$
are closed and the set $S_0=\Gamma_0$ consists of all the characters of $G=SL(2,\R)$. In addition, as defined in Section \ref{rep}, for every $i \in \{0,...,8\}$, there exists one common Hilbert space $\H_i$ that all the representations in $\Gamma_i$ act on. \\
Therefore, Condition 1 of Definition \ref{NCDL} is fulfilled. \\
~\\
As every semisimple Lie group with a finite center meets the CCR-condition, Condition 2 of Definition \ref{NCDL} is fulfilled as well. Thus, Condition 3 remains to be shown:

\section{Condition 3(a)}\label{cond 3a}
For the proof of Condition 3(a), as well as for the proof of Condition 3(b), some preliminaries are needed: 

\begin{lemma} \label{density-a}
~\\
Let $i \in \{0,...,8\}$, let $M$ be a dense subset of $C^*(G)$, let $\tilde{\nu}:CB(S_{i-1}) \to \B(\H_i)$ be a linear map bounded by $c \| \cdot \|_{S_{i-1}}$ for a constant $c>0$ and let $(\pi_j)_{j \in \N}$ be representations in $\Gamma_i$ such that
$$\Big\|\pi_j(h)-\tilde{\nu} \big(\F(h)\res{S_{i-1}} \big) \Big\|_{op} \overset{j \to \iy}{\longrightarrow}0~~~\forall h \in M.$$
Then
$$ \Big\|\pi_j(a)-\tilde{\nu} \big(\F(a)\res{S_{i-1}} \big) \Big\|_{op} \overset{j \to \iy}{\longrightarrow}0~~~\forall a \in C^*(G).$$
\end{lemma}
~\\
Proof: \\
Let $a \in C^*(G)$ and $\varepsilon >0$. \\
Since $M$ is dense in $C^*(G)$, there exists $h \in M$ such that $\|h-a\|_{C^*(G)} < \frac{\varepsilon}{3 \tilde{c}}$ for $\tilde{c}:=\max\{1,c\}$. \\
Furthermore, there exists $J( \varepsilon)>0$ in such a way that for all integers $j \geq J( \varepsilon)$ one has $\Big\|\pi_j(h)-\tilde{\nu} \big(\F(h)\res{S_{i-1}} \big) \Big\|_{op} < \frac{\varepsilon}{3}$. \\
Then, for all $j \geq J( \varepsilon)$
\begin{eqnarray*}
\Big\|\pi_j(a)-\tilde{\nu} \big(\F(a)\res{S_{i-1}} \big) \Big\|_{op} &\leq& \|\pi_j(a)-\pi_j(h)\|_{op}~+~ \Big\|\pi_j(h)-\tilde{\nu} \big(\F(h)\res{S_{i-1}} \big) \Big\|_{op}\\
&&+~ \Big\|\tilde{\nu} \big(\F(h)\res{S_{i-1}} \big)-\tilde{\nu} \big(\F(a)\res{S_{i-1}} \big) \Big\|_{op}.
\end{eqnarray*}
By assumption, $\Big\|\pi_j(h)-\tilde{\nu} \big(\F(h)\res{S_{i-1}} \big) \Big\|_{op} < \frac{\varepsilon}{3}$. In addition, as $\pi_j$ is a homomorphism
$$\|\pi_j(a)-\pi_j(h)\|_{op}~=~\|\pi_j(a-h)\|_{op}~\leq~\sup \limits_{\tilde{\pi} \in \wh{G}} \big\| \tilde{\pi}(a-h) \big\|_{op}~=~\|a-h\|_{C^*(G)}<\frac{\varepsilon}{3\tilde{c}}.$$
Moreover, 
\begin{eqnarray*}
\Big\|\tilde{\nu} \big(\F(h)\res{S_{i-1}} \big)-\tilde{\nu} \big(\F(a)\res{S_{i-1}} \big) \Big\|_{op}
&=& \Big\|\tilde{\nu} \big(\F(h)\res{S_{i-1}}-\F(a)\res{S_{i-1}} \big) \Big\|_{op} \\
&\leq&c~ \big\| \F(h)\res{S_{i-1}}-\F(a)\res{S_{i-1}} \big\|_{S_{i-1}} \\
&=&c ~\big\| \F(h-a)\res{S_{i-1}} \big\|_{S_{i-1}} \\
&=& c~\sup \limits_{\tilde{\pi} \in S_{i-1}} \big\| \F(h-a) (\tilde{\pi}) \big\|_{op} \\
&\leq& c~\sup \limits_{\tilde{\pi} \in \wh{G}} \big\| \tilde{\pi}(h-a) \big\|_{op} \\
&=&c ~\big\| h-a \big\|_{C^*(G)}< \frac{\varepsilon}{3}.
\end{eqnarray*}
Thus,
$$\Big\|\pi_j(a)-\tilde{\nu} \big(\F(a)\res{S_{i-1}} \big) \Big\|_{op} ~<~\varepsilon$$
and the claim is shown. \\
\qed
~\\
~\\

Define the $K \times K$-representation $\pi_{K \times K}$ on the space $V_{\pi_{K \times K}}:=C_0^{\iy}(G)$ of compactly supported $C^{\iy}(G)$-functions as 
\begin{eqnarray*}
\pi_{K \times K}&:& K \times K \to \B \big(C_0^{\iy}(G) \big),\\
&&\pi_{K \times K}(k_1,k_2)h(g):=h(k_1^{-1}g k_2)~~~\forall (k_1,k_2) \in K \times K~\forall h \in C_0^{\iy}(G)~\forall g \in G.
\end{eqnarray*}
For all $l,n \in \Z$ 
\begin{eqnarray*}
V_{\pi_{K \times K}}(l,n)&=& \Big\{h \in C_0^{\iy}(G) \big|~\pi_{K \times K} \big(k_{\va_1},k_{\va_2} \big)h= e^{i l \va_1 +i n\va_2}h~\forall \va_1,\va_2 \in [0,2 \pi) \Big\} \\
&=& \Big\{h \in C_0^{\iy}(G)\big|~ h \big(k_{\va_1}^{-1}g k_{\va_2} \big)= e^{i l \va_1 +i n\va_2}h(g)~\forall \va_1,\va_2 \in [0,2 \pi) \Big\}.
\end{eqnarray*}
Then the algebraic direct sum $\bigoplus \limits_{l,n \in \Z} V_{\pi_{K \times K}}(l,n)$ is dense in $V_{\pi_{K \times K}}=C_0^{\iy}(G)$ with respect to the $L^1(G)$-norm and as $\| \cdot \|_{C^*(G)} \leq \| \cdot \|_{L^1(G)}$ on $L^1(G)$, it is dense with respect to the $C^*(G)$-norm as well. $C_0^{\iy}(G)$ in turn is dense in $C^*(G)$. Hence, the algebraic direct sum $\bigoplus \limits_{l,n \in \Z} V_{\pi_{K \times K}}(l,n)$ is also dense in $C^*(G)$. \\
Let $p_{l,n}$ the projection going from $V_{\pi_{K \times K}}$ to $V_{\pi_{K \times K}}(l,n)$ defined in the following way: \\
For $h \in V_{\pi_{K \times K}}=C_0^{\iy}(G)$ and $g \in G$
$$p_{l,n}(h)(g)~:=~ \frac{1}{|K|^2} \int \limits_{K \times K} h \big(k_{\va_1} g k_{\va_2}^{-1} \big) e^{i l \va_1} e^{i n \va_2} d \big( k_{\va_1},k_{\va_2} \big).$$
Similar as shown in Lemma \ref{density-a} above, in order to prove that a sequence of representations applied to general elements $a \in C^*(G)$ converges to a representation applied to $a$, it suffices to show this convergence for elements $h$ in a dense subset $M \subset C^*(G)$. Therefore and with Lemma \ref{density-a}, due to the density discussed above, instead of dealing with a general $a \in C^*(G)$, the calculations in Sections \ref{cond 3a} and \ref{cond 3b} can be accomplished with a function $h \in C_0^{\iy}(G)$ fulfilling $h=p_{l,-n}(h)$ for some integers $l,n \in \Z$.

\begin{lemma} \label{rel 3(b)(iii)}
~\\
For $f \in L^2(K)_+$ and $h \in C_0^{\iy}(G)$ with $h=p_{l,-n}(h)$ for $l, n \in \Z$, one gets
\begin{eqnarray*}
\P^{+,u}(h)(f)~=~\P^{+,u}(h) \big(p_{n}(f) \big)~=~p_{l} \Big(\P^{+,u}(h) \big(p_{n}(f) \big) \Big)~~~\forall u \in \C.
\end{eqnarray*}
Furthermore, if $h \not =0$ the integers $l$ and $n$ must be even.
\end{lemma}

This lemma is a standard formula. Its proof consists of a simple calculation and will thus be skipped. \\
~\\
Now, the proof of Condition 3(a) can be executed. \\
Condition 3(a) is obvious for the sets $\Gamma_i$ for $i \in \{0,...,5\}$, as these are discrete sets. \\
~\\
So, let $v_j,v \in [0, \iy)$ and $\big(\P^{+,iv_j}\big)_{j \in \N}$ a sequence in $\Gamma_6$ converging to $\P^{+,iv}$. Then $v_j \overset{j \to \iy}{\longrightarrow}v$. Hence, let $h \in C_0^{\iy}(G)$ be supported in the compact set $\tilde{K} \in G$. Then, very similar as in the proof of (\ref{hilfskonvergenz}), there is a compact set $I \subset \R$ such that
\begin{eqnarray} \label{op-norm-h}
\nn &&\big\| \P^{+,iv_j}(h)- \P^{+,iv}(h)\big\|_{op}^2\\
\nn &=& \sup \limits_{\substack{f \in L^2(K)_+\\ \|f\|_2=1}}~ \int \limits_K \bigg| \int \limits_{\tilde{K}} h(g) \Big( \P^{+,iv_j}(g)(f)(k)- \P^{+,iv}(g)(f)(k)\Big)dg \bigg|^2dk\\
\nn &\leq& \sup \limits_{x \in I} \Big|e^{-iv_j x}-e^{-iv x}\Big|^2 \sup \limits_{\substack{f \in L^2(K)_+\\ \|f\|_2=1}}~\int \limits_K \bigg( \int \limits_{\tilde{K}} |h(g)|~ e^{- \rho H(g^{-1}k)} \Big|f \Big(\kappa \big(g^{-1}k \big)\Big)\Big| dg \bigg)^2dk \\
\nn &\overset{\text{H\"older}}{\leq}& \sup \limits_{x \in I} \Big|e^{-iv_j x}-e^{-iv x}\Big|^2 ~\|h\|_{L^2(G)}^2 \sup \limits_{\substack{f \in L^2(K)_+ \\ \|f\|_2=1}}~\int \limits_K \int \limits_{\tilde{K}} e^{-2 \rho H(g^{-1}k)}\Big|f \Big(\kappa \big(g^{-1}k \big)\Big)\Big|^2 dg dk \\
&\overset{Lemma~\ref{intabsch}}{=}&\big| \tilde{K} \big| \sup \limits_{x \in I} \Big|e^{-iv_j x}-e^{-iv x}\Big|^2 ~\|h\|_{L^2(G)}^2
\overset{j \to \iy}{\longrightarrow}0,
\end{eqnarray}
as $v_j \overset{j \to \iy}{\longrightarrow}v$.\\
Because of the density of $C_0^{\iy}(G)$ in $C^*(G)$, one gets the desired convergence for $a \in C^*(G)$. \\
The reasoning is the same for $\Gamma_{7}$. \\
~\\ 
For $\Gamma_{8}$, let $u_j,u \in (0,1)$ and $\big({\cal C}^{u_j}\big)_{j \in \N}$ a sequence in $\Gamma_{8}$ converging to ${\cal C}^{u} \in \Gamma_{8}$. Then $u_j \overset{j \to \iy}{\longrightarrow}u$.
Moreover, let $h \in C^*(G)$ and let $l, n \in \Z$ such that $h=p_{l,-n}(h)$, as discussed at the beginning of this section. \\
Let $f \in L^2(K)_+$ with $\| f \|_{L^2(K)}=1$. \\
Since $K_{\tilde{u}}^{-1}$ commutes with $p_n$, by Lemma \ref{rel 3(b)(iii)} one has for every $\tilde{u} \in (0,1)$
\begin{eqnarray}\label{cu-rechnung}
\nn {\cal C}^{\tilde{u}}(h)(f)
&=& K_{\tilde{u}} \circ \P^{+,\tilde{u}}(h) \circ K_{\tilde{u}}^{-1}(f) 
~=~ K_{\tilde{u}} \circ p_l \bigg(  \P^{+,\tilde{u}}(h) \Big( p_n \big( K_{\tilde{u}}^{-1}(f) \big) \Big) \bigg)\\
\nn &=& d_l \big( \tilde{u} \big)~ p_l \bigg(  \P^{+,\tilde{u}}(h) \Big( K_{\tilde{u}}^{-1} \big( p_n(f) \big) \Big) \bigg)
~=~ \frac{d_l \big( \tilde{u} \big)}{d_n \big( \tilde{u} \big)}~ p_l \Big( \P^{+,\tilde{u}}(h) \big( p_n(f) \big) \Big) \\
&=& \sqrt{\frac{c_{l} \big( \tilde{u})}{c_{n} \big(\tilde{u} \big)}} ~ \P^{+,\tilde{u}}(h)(f).
\end{eqnarray}
Hence, 
\begin{eqnarray} \label{3afuergamma8}
\big\|{\cal C}^{u_j}(h)- {\cal C}^{u}(h) \big\|_{op}^2
&=& \Bigg\| \sqrt{\frac{c_{l}(u_j)}{c_{n}(u_j)}} ~ \P^{+,u_j}(h) - 
\sqrt{\frac{c_l(u)}{c_{n}(u)}} ~\P^{+,u}(h) \Bigg\|_{op}^2 \overset{j \to \iy}{\longrightarrow}0,
\end{eqnarray}
since $u_j \overset{j \to \iy}{\longrightarrow} u$ and with the same reasoning as in (\ref{op-norm-h}).\\
Because of the density of $C_0^{\iy}(G)$ in $C^*(G)$, one gets the desired convergence for $a \in C^*(G)$ and thus, the claim is also shown for $\Gamma_{8}$. 

\section{Condition 3(b)} \label{cond 3b}
Now, only Condition 3(b) remains to be shown. This is the most complicated part of the proof of the different conditions listed in Definition \ref{NCDL}.\\
~\\
The setting of a sequence $(\gamma_j)_{j \in \N} $ in $\Gamma_i$ converging to a limit set contained in 
$S_{i-1}=\bigcup \limits_{l<i}\Gamma_l$ regarded in Condition 3(b) can only occur in the following cases: 
\begin{enumerate}[(i)]
\item $(\gamma_j)_{j \in \N}=\big(\P^{-,iv_j}\big)_{j \in \N}$ is a sequence in $\Gamma_{7}$ whose limit set is $\Gamma_3 \cup \Gamma_4=\big\{\D_+, \D_-\big\}$.
\item $(\gamma_j)_{j \in \N}=\big({\cal C}^{u_j}\big)_{j \in \N}$ is a sequence in $\Gamma_{8}$ whose limit set is $\big\{\P^{+,0} \big\} \subset \Gamma_{7}$.
\item $(\gamma_j)_{j \in \N}=\big({\cal C}^{u_j}\big)_{j \in \N}$ is a sequence in $\Gamma_{8}$ whose limit set is $\Gamma_0 \cup \Gamma_1 \cup \Gamma_2=\big\{\F_1,\D_1^+, \D_1^-\big\}$.
\end{enumerate}

For $i \in \{0,...,6\}$ the sets $\Gamma_i$ are closed and thus the regarded situation cannot appear for sequences $(\gamma_j)_{j \in \N} $ in $\Gamma_i$ for $i \in \{0,...,6\}$. \\
~\\
Since all the sequences regarded in the cases (i), (ii) and (iii) are properly converging, the transition to subsequences will be omitted. \\
~\\
~\\
First, regard Case (i):\\
Let $\big(\P^{-,iv_j}\big)_{j \in \N}$ be a sequence in $\Gamma_{7}$ whose limit set is $\big\{\D_+, \D_-\big\}$. 
As $\P^{-,iv_j} \overset{j \to \iy}{\longrightarrow} \big\{\D_+, \D_-\big\}$, it follows that $v_j \overset{j \to \iy}{\longrightarrow}0$. \\
Now, a bounded, linear and involutive mapping
$$\tilde{\nu}_j:CB(S_6) \to \B\big(L^2(K)_-\big)$$
fulfilling
$$\lim \limits_{j \to \iy} \Big\| \P^{-, i v_j}(a)-\tilde{\nu}_j \big(\F(a)\res{S_6} \big) \Big\|_{op}=0~~~\forall a \in C^*(G)$$
has to be defined. \\
Since in this and in the following cases this mapping will not depend on $j$, it will from now on be denoted by $\tilde{\nu}$ instead of $\tilde{\nu}_j$. \\
Let $p_+$ be the projection from $L^2(K)_{-}$ to the space $\H_{\D_+}$ and $p_-$ the projection from $L^2(K)_{-}$ to $\H_{\D_-}$. Then, by Remark \ref{rem1} one has $L^2(K)_-=\H_{\D_+} \oplus \H_{\D_-}$, i.e. $\id \res{L^2(K)_-}=p_+ + p_-$. \\
Now, let
$$\tilde{\nu}(\psi):=\tilde{\nu}_{\{\D_+, \D_-\}}(\psi):= \psi(\D_+) \circ p_+ + \psi(\D_-) \circ p_-~~~\forall \psi \in CB(S_6).$$
This is well-defined, as $\D_+, \D_- \in S_6$, and furthermore $\tilde{\nu}(\psi) \in \B \big(L^2(K)_-\big)$. \\
The linearity of the mapping $\tilde{\nu}$ is clear. For the involutivity, let $\psi \in CB(S_6)$. Then, since $p_+$ and $p_-$ equal the identity on the image of $\psi(\D_+)$ and respectively $\psi(\D_-)$,
\begin{eqnarray*}
\big(\tilde{\nu}(\psi)\big)^*&=&\Big( \psi(\D_+)\circ p_+\Big)^*+\Big( \psi(\D_-)\circ p_-\Big)^*~=~\Big(p_+ \circ \psi(\D_+)\circ p_+\Big)^*+\Big(p_- \circ \psi(\D_-)\circ p_-\Big)^* \\
&=& p_+^* \circ \psi^*(\D_+)\circ p_+^* + p_-^* \circ \psi^*(\D_-)\circ p_-^* ~=~ p_+ \circ \psi^*(\D_+)\circ p_+ + p_- \circ \psi^*(\D_-)\circ p_- \\
&=& \psi^*(\D_+)\circ p_+ + \psi^*(\D_-)\circ p_-~=~\tilde{\nu}(\psi^*).
\end{eqnarray*}
To show that $\tilde{\nu}$ is bounded, again let $\psi \in CB(S_6)$:
\begin{eqnarray*}
\big\|\tilde{\nu}(\psi) \big\|_{op}&=&\big\| \psi(\D_+)\circ p_+ + \psi(\D_-)\circ p_- \big\|_{op}
~=~\max \Big\{ \big\| \psi(\D_+)\big\|_{op}, \big\| \psi(\D_-)\big\|_{op} \Big\} \\
&\leq& \sup \limits_{\gamma \in S_6} \| \psi(\gamma)\|_{op}~=~\| \psi \|_{S_6}. 
\end{eqnarray*}
Now, only the demanded convergence remains to be shown: \\
For $h \in C_0^{\iy}(G)$, one has
\begin{eqnarray*}
\Big\| \P^{-, i v_j}(h)-\tilde{\nu} \big(\F(h) \res{S_6} \big) \Big\|^2_{op}
&=& \bigg\| \P^{-,iv_j}(h)-\Big( \F(h)(\D_+)\circ p_++ \F(h)(\D_-)\circ p_- \Big) \bigg\|_{op}^2\\
&=&\bigg\| \P^{-,iv_j}(h)-\Big( \D_+(h)\circ p_++\D_-(h)\circ p_- \Big) \bigg\|_{op}^2\\
&=&\bigg\| \P^{-,iv_j}(h)-\Big( \P^{-,0}(h)\circ p_++\P^{-,0}(h)\circ p_- \Big) \bigg\|_{op}^2\\
&=& \big\| \P^{-,iv_j}(h)-\P^{-,0}(h) \big\|_{op}^2 \overset{j \to \iy}{\longrightarrow}0
\end{eqnarray*}
as in (\ref{op-norm-h}) and since $v_j \overset{j \to \iy}{\longrightarrow}0$. \\
Again, because of the density of $C_0^{\iy}(G)$ in $C^*(G)$ and with Lemma \ref{density-a}, one gets the desired convergence for $a \in C^*(G)$. \\
\newpage
Now, regard Case (ii):\\
Let $\big({\cal C}^{u_j}\big)_{j \in \N}$ be a sequence in $\Gamma_{8}$ whose limit set is $\big\{\P^{+,0}\big\}$. Thus, $u_j \overset{j \to \iy}{\longrightarrow}0$. \\
Here, a bounded, linear and involutive mapping
$$\tilde{\nu}:CB(S_{7}) \to \B \big(L^2(K)_+ \big)$$
fulfilling
$$\lim \limits_{j \to \iy} \Big\| {\cal C}^{u_j}(a)-\tilde{\nu} \big(\F(a) \res{S_{7}} \big) \Big\|_{op}=0~~~\forall a \in C^*(G)$$
is needed: \\
Define
$$\tilde{\nu}(\psi):=\tilde{\nu}_{\P^{+,0}}(\psi):=\psi \big(\P^{+,0} \big)~~~\forall \psi \in CB(S_{7}).$$
$\tilde{\nu}(\psi) \in \B \big(L^2(K)_+\big)$ for every $\psi \in CB(S_{7})$ and $\tilde{\nu}$ is well-defined, as $\P^{+,0} \in S_{7}$.  \\
The linearity and the involutivity of $\tilde{\nu}$ are clear. \\
For the boundedness of $\tilde{\nu}$, let $\psi \in CB(S_{7})$. Then 
\begin{eqnarray*}
\big\|\tilde{\nu}(\psi) \big\|_{op} ~=~ \big\| \psi \big(\P^{+,0} \big) \big\|_{op}
~\leq~  ~\sup \limits_{\gamma \in S_{7}} \| \psi(\gamma)\|_{op}~=~\| \psi \|_{S_{7}}. 
\end{eqnarray*}
Again, it remains to show the demanded convergence: \\
So, let $h \in C_0^{\iy}(G)$. Then one can assume again that there exist $l, n \in \Z$ such that $h=p_{l,-n}(h)$. Since $\lim \limits_{u \to 0} \frac{c_{\tilde{n}}(u)}{c_{n'}(u)}=1$ for all $\tilde{n},n' \in \Z$ by (\ref{cn}) and Lemma \ref{id}, one gets with (\ref{cu-rechnung})
\begin{eqnarray*}
\Big\| {\cal C}^{u_j}(h)~-~\tilde{\nu} \big(\F(h) \res{S_{7}} \big) \Big\|_{op}^2
&=& \Bigg\| \sqrt{\frac{c_l(u_j)}{c_n(u_j)}}~ \P^{+,u_j}(h) ~-~ \F(h) \big(\P^{+,0} \big) \Bigg\|_{op}^2 \\
&=& \Bigg\| \sqrt{\frac{c_l(u_j)}{c_n(u_j)}}~ \P^{+,u_j}(h) ~-~  \P^{+,0}(h) \Bigg\|_{op}^2 \overset{j \to \iy}{\longrightarrow}0,
\end{eqnarray*}
since $u_j \overset{j \to \iy}{\longrightarrow}0$ and with the same arguments as in the section above. \\
The desired convergence for $a \in C^*(G)$ follows. \\
~\\
Last, regard Case (iii):\\
Let $\big({\cal C}^{u_j}\big)_{j \in \N}$ be a sequence in $\Gamma_{8}$ whose limit set is $\big\{\F_1,\D_1^+, \D_1^-\big\}$. This means that $u_j \overset{j \to \iy}{\longrightarrow}1$. \\
Again, a bounded, linear and involutive mapping
$$\tilde{\nu}:CB(S_{7}) \to \B \big(L^2(K)_+\big)$$
fulfilling
$$\lim \limits_{j \to \iy} \Big\| {\cal C}^{u_j}(a)-\tilde{\nu} \big(\F(a) \res{S_{7}} \big) \Big\|_{op}=0~~~\forall a \in C^*(G)$$
is needed: \\
For this, let $p_+$ the projection from $L^2(K)_{+}$ to the space $\big\{f \in L^2(K)_+ |~p_n(f)=0 ~\forall n \leq 0 \big\}$ and $p_-$ the projection from $L^2(K)_{+}$ to $\big\{f \in L^2(K)_+ |~p_n(f)=0 ~\forall n \geq 0 \big\}$. Then, since
\begin{eqnarray*}
L^2(K)_+&=&\big\{f \in L^2(K)_+ |~p_n(f)=0 ~\forall n \leq 0 \big\}~+~\big\{f \in L^2(K)_+ |~p_n(f)=0 ~\forall n \geq 0 \big\} \\
&&+~\big\{f \in L^2(K)_+ |~p_n(f)=0 ~\forall n \not= 0 \big\},
\end{eqnarray*}
every $f \in L^2(K)_{+}$ can be written as $f=p_+(f)+p_-(f)+p_0(f)$.\\
Furthermore, let
$$d_{n,2}(1):= \lim \limits_{u \to 1} \sqrt{\frac{c_n(u)}{c_2(u)}}~~~\text{for all even}~n>0~~~\text{and}$$
$$d_{n,-2}(1):= \lim \limits_{u \to 1} \sqrt{\frac{c_n(u)}{c_{-2}(u)}}~~~\text{for all even}~n<0.$$ 
The existence of these limits follows with Lemma \ref{existence}(c) in Section \ref{Ku} and the analog statement for $\tilde{J}_{[1]}$. \\
Now, define the operators
$$ K_{(1)}: \H_{(1)} \to L^2(K)_+~~~\text{by}~~~{K_{(1)}}\res{\H_{\P^{+,1}}(n)}:=d_{n,2}(1) \cdot \id\res{  \H_{\P^{+,1}}(n)}~~~\text{for all even}~n>0$$
and
$$ K_{[1]}: \H_{[1]} \to L^2(K)_+~~~~\text{by}~~~{K_{[1]}}\res{ \H_{\P^{+,1}}(n)}:=d_{n,-2}(1) \cdot \id\res{ \H_{\P^{+,1}}(n)}~~~\text{for all even}~n<0.$$
By definition, these operators are linear and they are unitary by the construction of the scalar products defined on $\H_{(1)}$ and $\H_{[1]}$. Moreover, like $\tilde{J}_{(1)}$ and $\tilde{J}_{(-1)}$, they are injective
and, as proved in Remark \ref{commute} for the operators $J_u$, $K_{(1)}$ and $K_{[1]}$ commute with the projections $p_n$ for all $n \in \N$. \\
Moreover, one can easily see that
$$K_{(1)}\big( \H_{(1)}\big) ~=~ \big\{f \in L^2(K)_+ |~p_n(f)=0 ~\forall n \leq 0 \big\}~=~p_+ \big(L^2(K)_+ \big)~~~\text{and}$$
$$K_{[1]}\big( \H_{[1]}\big) ~=~ \big\{f \in L^2(K)_+ |~p_n(f)=0 ~\forall n \geq 0 \big\}~=~p_- \big(L^2(K)_+ \big).$$
Hence, one can build the inverse of the operators $K_{(1)}$ and $K_{[1]}$ on the image of $p_+$ or respectively $p_-$. Therefore, one can define
\begin{eqnarray*}
\tilde{\nu}(\psi)~:=~\tilde{\nu}_{\{\F_1,\D_1^+, \D_1^-\}}(\psi)&:= &K_{(1)} \circ \psi(\D_1^+) \circ K_{(1)}^{-1}\circ p_+ ~+~ K_{[1]} \circ \psi(\D_1^-)\circ K_{[1]}^{-1} \circ p_- \\
&&+~ \psi(\F_1) \circ p_0 ~~~\forall \psi \in CB(S_{7}).
\end{eqnarray*}
The mapping $\tilde{\nu}$ is well-defined, since $\D_1^+, \D_1^-,\F_1 \in S_{7}$, and $\tilde{\nu}(\psi) \in \B \big(L^2(K)_+\big)$ for every $\psi \in CB(S_{7})$. In addition, its linearity is clear again. \\
For the involutivity and the boundedness, let $\psi \in CB(S_{7})$. Then, as $K_{(1)}$ and $K_{[1]}$ are unitary and as $p_+$, $p_-$ and $p_0$ equal the identity on the image of $K_{(1)}$, $K_{[1]}$ and respectively $\psi(\F_1)$, 
\begin{eqnarray*}
&&\big(\tilde{\nu}(\psi) \big)^* \\
&=& \Big( K_{(1)} \circ \psi(\D_1^+) \circ K_{(1)}^{-1}\circ p_+~ +~K_{[1]} \circ \psi(\D_1^-)\circ K_{[1]}^{-1} \circ p_- ~+~ \psi(\F_1) \circ p_0 \Big)^* \\
&=& \Big( p_+ \circ K_{(1)} \circ \psi(\D_1^+) \circ K_{(1)}^{-1}\circ p_+~ +~ p_- \circ K_{[1]} \circ \psi(\D_1^-)\circ K_{[1]}^{-1} \circ p_- ~+~p_0 \circ \psi(\F_1) \circ p_0 \Big)^* \\
&=&  p_+ \circ K_{(1)} \circ \psi^*(\D_1^+) \circ K_{(1)}^{-1}\circ p_+~ +~p_- \circ K_{[1]} \circ \psi^*(\D_1^-)\circ K_{[1]}^{-1} \circ p_- ~+~p_0 \circ \psi^*(\F_1) \circ p_0 \\
&=& \tilde{\nu} \big(\psi^* \big).
\end{eqnarray*}
Furthermore, since $\big\| K_{(1)} \big\|_{op} \big\| K_{(1)} ^{-1} \big\|_{op}= \big\| K_{[1]} \big\|_{op} \big\| K_{[1]}^{-1} \big\|_{op} = 1$, one gets
\begin{eqnarray*}
\big\|\tilde{\nu}(\psi) \big\|_{op}
&=& \Big\|K_{(1)} \circ \psi(\D_1^+) \circ K_{(1)}^{-1}\circ p_+~ +~K_{[1]} \circ \psi(\D_1^-)\circ K_{[1]}^{-1} \circ p_- ~+~\psi(\F_1) \circ p_0 \Big\|_{op} \\
&=& \max \bigg\{ \Big\| K_{(1)} \circ \psi(\D_1^+) \circ K_{(1)}^{-1}\Big\|_{op},~ \Big\| K_{[1]} \circ \psi(\D_1^-)\circ K_{[1]}^{-1}\Big\|_{op},~
 \big\| \psi(\F_1) \big\|_{op} \bigg\} \\
&=& \max \Big\{ \big\| \psi(\D_1^+) \big\|_{op},~ \big\| \psi(\D_1^-) \big\|_{op},~
 \big\| \psi(\F_1) \big\|_{op} \Big\} \\
&\leq& ~\sup \limits_{\gamma \in S_{7}} \| \psi(\gamma)\|_{op}~=~\| \psi \|_{S_{7}}. 
\end{eqnarray*}
For the demanded convergence, let $h \in C_0^{\iy}(G)$. As above in the proof of (ii) and Condition 3(a), one can assume that there exist $l, n \in \Z$ such that $h=p_{l,-n}(h)$. \\
Let $f \in L^2(K)_+$ with $\| f \|_{L^2(K)}=1$. \\
Since $K_{(1)}^{-1}$ and $K_{[1]}^{-1}$ commute with $p_n$, similar as in the proof of (\ref{cu-rechnung}), by Lemma \ref{rel 3(b)(iii)} one gets
\begin{eqnarray*}
&& \tilde{\nu} \big(\F(h)\res{S_{7}} \big)(f) \\
&=& K_{(1)} \circ \F(h)(\D_1^+) \circ K_{(1)}^{-1} \circ p_+(f) ~+~ K_{[1]} \circ \F(h)(\D_1^-)\circ K_{[1]}^{-1} \circ p_-(f)~+ ~\F(h)(\F_1) \circ p_0(f)\\
&=& K_{(1)} \circ \D_1^+(h) \circ K_{(1)}^{-1} \circ p_+ (f)~+~ K_{[1]} \circ \D_1^-(h)\circ K_{[1]}^{-1} \circ p_-(f)~+ ~\F_1(h) \circ p_0 (f) \\
&=& K_{(1)} \circ \P^{+,1}(h) \circ K_{(1)}^{-1} \circ p_+ (f)~+~ K_{[1]} \circ \P^{+,1}(h) \circ K_{[1]}^{-1} \circ p_-(f)~+~ \P^{+,-1}(h) \circ p_0 (f) \\ 
&=& K_{(1)} \circ p_l \Bigg( \P^{+,1}(h) \bigg(p_n \Big( K_{(1)}^{-1} \circ p_+ (f) \Big)\bigg)\Bigg)~+~ K_{[1]} \circ p_l \Bigg( \P^{+,1}(h) \bigg(p_n \Big( K_{[1]}^{-1} \circ p_- (f) \Big)\bigg)\Bigg)\\
&&+ ~ p_l \bigg( \P^{+,-1}(h) \Big(p_n \big(p_0 (f) \big)\Big)\bigg)\\
&=& d_{l,2}(1)~ p_l \Bigg( \P^{+,1}(h) \bigg(K_{(1)}^{-1} \Big(p_n  \circ p_+ (f) \Big)\bigg)\Bigg)~+~ d_{l,-2}(1)~ p_l \Bigg( \P^{+,1}(h) \bigg(K_{[1]}^{-1} \Big(p_n  \circ p_- (f) \Big)\bigg)\Bigg) \\
&&+ ~p_l \bigg( \P^{+,-1}(h) \Big(p_n \circ p_0 (f)\Big)\bigg).
\end{eqnarray*}
Now, there are three cases to consider: In the first case $n >0$, in the second case $n<0$ and in the third case $n=0$. \\
So, first let $n>0$. Then
\begin{eqnarray*}
\tilde{\nu} \big(\F(h)\res{S_{7}} \big)(f)
&=& d_{l,2}(1)~ p_l \Bigg( \P^{+,1}(h) \bigg(K_{(1)}^{-1} \Big(p_n  \circ p_+ (f) \Big)\bigg)\Bigg) \\
&=& \frac{d_{l,2}(1)}{d_{n,2}(1)}~ p_l \Big( \P^{+,1}(h) \big(p_n (f) \big) \Big) 
~=~\lim \limits_{u \to 1} \sqrt{\frac{c_l(u)}{c_n(u)}}~ \P^{+,1}(h)(f).
\end{eqnarray*}
Therefore, joining this result with (\ref{cu-rechnung}), one gets
\begin{eqnarray} \label{3teile}
\nn \Big\| {\cal C}^{u_j}(h)~-~ \tilde{\nu} \big(\F(h)\res{S_{7}} \big) \Big\|_{op}^2
\nn &=& \Bigg\| \sqrt{\frac{c_l(u_j)}{c_n(u_j)}}~ \P^{+,u_j}(h) ~-~
\lim \limits_{u \to 1} \sqrt{\frac{c_l(u)}{c_n(u)}}~ \P^{+,1}(h) \Bigg\|_{op}^2 \overset{j \to \iy}{\longrightarrow}0,
\end{eqnarray}
since $u_j \overset{j \to \iy}{\longrightarrow} 1$ and with the same reasoning as in (\ref{op-norm-h}). \\
~\\
The proof for the case $n<0$ is the same as the one for the first case, hence it only remains to regard the case $n=0$: \\
Since by (\ref{c_n(1)})
$$\frac{c_l(1)}{c_0(1)}~=~\begin{cases} 1~~~\text{if}~l=0 \\ 0~~~\text{if}~l \not= 0, \end{cases}$$
one gets for $l \not=0$ with (\ref{cu-rechnung})
\begin{eqnarray*}
\Big\| {\cal C}^{u_j}(h)~-~ \tilde{\nu} \big(\F(h)\res{S_{7}} \big) \Big\|_{op}^2
&=& \Bigg\| \sqrt{\frac{c_l(u_j)}{c_0(u_j)}}~ \P^{+,u_j}(h) ~-~0 \Bigg\|_{op}^2 \overset{j \to \iy}{\longrightarrow}0,
\end{eqnarray*}
since $u_j \overset{j \to \iy}{\longrightarrow} 1$. \\
So, let $l=0$ and define $C_h:=\int \limits_G h(g)dg$. \\
Then $\F_1(h)=C_h \cdot \id$.\\
Furthermore, with Lemma \ref{intabsch}, one gets for $g \in G$
\begin{eqnarray*}
p_0 \circ \P^{+,1}(g) \circ p_0(f)&=& \frac{1}{|K|} \int \limits_K  \P^{+,1}(g) p_0(f)(k) dk \\
&=& \frac{p_0(f)}{|K|} \int \limits_K e^{-2 \rho H(g^{-1}k)} dk \\
&=& \frac{p_0(f)}{|K|} \int \limits_K e^{-2 \rho H(g^{-1}k)} \Big|1 \Big(\kappa\big(g^{-1}k \big) \Big)\Big|^2 dk\\
&=& \frac{p_0(f)}{|K|}~ \|1\|_{L^2(K)}~=~p_0(f).
\end{eqnarray*}
Thus, for $h$ one has
$$p_0 \circ \P^{+,1}(h) \circ p_0(f)~=~ \int \limits_G h(g) p_0(f) dg~=~C_h p_0(f)~=~\F_1(h) \circ p_0 (f).$$
Therefore, for $l=n=0$, by (\ref{cu-rechnung})
\begin{eqnarray*}
\Big\| {\cal C}^{u_j}(h)~-~ \tilde{\nu} \big(\F(h)\res{S_{7}} \big) \Big\|_{op}^2
&=& \Bigg\| p_0 \circ \P^{+,u_j}(h) \circ p_0 ~-~p_0 \circ \P^{+,1}(h) \circ p_0 \Bigg\|_{op}^2 \overset{j \to \iy}{\longrightarrow}0,
\end{eqnarray*}
since $u_j \overset{j \to \iy}{\longrightarrow} 1$. \\
Hence the claim is also shown in the case $n=0$ and thus
$$\Big\| {\cal C}^{u_j}(h)~-~ \tilde{\nu} \big(\F(h)\res{S_{7}} \big) \Big\|_{op}^2\overset{j \to \iy}{\longrightarrow}0,$$
as demanded. \\
Again, the desired convergence for $a \in C^*(G)$ follows. \\
\qed

\section{Result}\label{result}
Having now verified all the conditions listed in Section \ref{cond}, Theorem \ref{NCDL - sl2r} is proved and the $C^*$-algebra of $G=SL(2, \R)$ can therefore be characterized as in Theorem \ref{theo1 - sl2r} with the sets $\Gamma_i$ and $S_i$ and the Hilbert spaces $\H_i$ for $i \in \{0,...,8\}$ defined in Section \ref{gammas} and Section \ref{rep} and the mappings $\tilde{\nu}$ constructed in Section \ref{cond 3b}. \\
~\\
Let for a topological Hausdorff space $V$ and a $C^*$-algebra $B$, $C_{\iy}(V,B)$ be the \mbox{$C^*$-algebra} of all continuous functions defined on $V$ with values in $B$ that are vanishing at infinity. Then, from Theorem \ref{theo1 - sl2r}, one can deduce more concretely the following result for $G=SL(2, \R)$:

\begin{theorem}\label{thm genauer - sl2r}
~\\
Let the operator $p_+$ be the projection from $L^2(K)_{\pm}$ to $\big\{f \in L^2(K)_{\pm} |~p_n(f)=0 ~\forall n \leq 0 \big\}$, $p_-$ the projection from $L^2(K)_{\pm}$ to the space $\big\{f \in L^2(K)_{\pm} |~p_n(f)=0 ~\forall n \geq 0 \big\}$ and $p_0$ the projection from $L^2(K)_{+}$ to $\big\{f \in L^2(K)_{+} |~p_n(f)=0 ~\forall n \not= 0 \big\}=\C$. \\
Then the $C^*$-algebra $C^*(G)$ of $G=SL(2, \R)$ is isomorphic to the direct sum of $C^*$-algebras
\begin{eqnarray*}
&&\bigg\{ F \in C_{\iy} \Big(i [0, \iy) \cup [0,1], \K \big(L^2(K)_+ \big) \Big) \big|~F(1)~\text{commutes with}~p_+,~p_-~\text{and}~p_0\bigg\} \\
&\oplus& \bigg\{ F \in C_{\iy} \Big(i [0, \iy), \K \big(L^2(K)_- \big) \Big) \big|~F(0)~\text{commutes with}~p_+~\text{and}~p_-\bigg\} \\
&\oplus& C_{\iy}\Big(\Z \setminus \{-1,0,1\}, \K \big(\H_{\D} \big) \Big)
\end{eqnarray*}
for the infinite-dimensional and separable Hilbert space $\H_{\D}$ fixed in Section \ref{rep}. 
\end{theorem}

Proof: \\
The unitary dual of $C^*(G)$ or respectively of $G$ is given by the disjoint union 
$$\wh{G}_{even} ~\overset{.}{\cup}~ \wh{G}_{odd}~ \overset{.}{\cup}~\wh{G}_{discrete},$$
where the set $\wh{G}_{even}$ consists of the even representations $\P^{+,iv}$ for $v \in [0, \iy)$, ${\cal C}^u$ for $u \in (0,1)$, $\D_1^+$, $\D_1^-$ and $\F_1$, the set $\wh{G}_{odd}$ consists of the odd representations $\P^{-,iv}$ for $v \in (0, \iy)$, $\D_+$ and $\D_-$ and the set $\wh{G}_{discrete}$ consists of the even or respectively odd representations $\D_m^+$ for $m \in \N_{>1}$ and $\D_m^-$ for $m \in \N_{>1}$. \\
These three listed sets of representations are topologically separated from each other (see Section \ref{top}). \\
Mapping $\P^{+,iv}$ for $v \in [0, \iy)$ to $iv$, ${\cal C}^u$ for $u \in (0,1)$ to $u$ and $\D_1^+$, $\D_1^-$ and $\F_1$ to $1$, one gets a surjection from $\wh{G}_{even}$ to the set $i[0, \iy) \cup [0,1]=:I_1$. \\
Furthermore, mapping $\P^{-,iv}$ for $v \in (0, \iy)$ to $iv$ and $\D_+$ and $\D_-$ to $0$, one gets a surjection from $\wh{G}_{odd}$ to $i[0, \iy)=:I_2$. \\
Last, mapping $\D_m^+$ for $m \in \N_{>1}$ to $m$ and $\D_m^-$ for $m \in \N_{>1}$ to $-m$, one gets a surjection from $\wh{G}_{discrete}$ to $\Z \setminus \{-1,0,1\}=:I_3$. \\
Hence, one regards the three sets
$$I_1~=~i[0, \iy) \cup [0,1], ~~~I_2~=~i [0, \iy)~~~\text{and}~~~I_3~=~\Z \setminus \{-1,0,1\}.$$
~\\
In order to prove this theorem, it has to be shown that for every operator field $\va= \F(a)$ for $a \in C^*(G)$ that fulfills the properties listed in Theorem \ref{theo1 - sl2r}, there exists a mapping \mbox{$F_a^1 \in \bigg\{ F \in C_{\iy} \Big(i [0, \iy) \cup [0,1], \K \big(L^2(K)_+ \big) \Big) \big|~F(1)~\text{commutes with}~p_+,~p_-~\text{and}~p_0\bigg\}=:P_1$}, a mapping $F_a^2 \in \bigg\{ F \in C_{\iy} \Big(i [0, \iy), \K \big(L^2(K)_- \big) \Big) \big|~F(0)~\text{commutes with}~p_+~\text{and}~p_-\bigg\}=:P_2$ and a mapping $F_a^3 \in C_{\iy}\Big(\Z \setminus \{-1,0,1\}, \K \big(\H_{\D} \big) \Big)=:P_3$. \\
On the other hand, for every $F_1 \in P_1$, every $F_2 \in P_2$ and every $F_3 \in P_3$ one has to construct an operator field $\va_{F_1,F_2,F_3}$ over $\wh{G}$ that fulfills the properties of Theorem \ref{theo1 - sl2r}. Since the above mentioned sets of representations $\wh{G}_{even}$, $\wh{G}_{odd}$ and $\wh{G}_{discrete}$ are topologically separated from each other, it suffices thus to define three different operator fields $\va_{F_1}$ over $\wh{G}_{even}$, $\va_{F_2}$ over $\wh{G}_{odd}$ and $\va_{F_3}$ over $\wh{G}_{discrete}$.\\
~\\
For every $a \in C^*(G)$ define a function $F_a^1:I_1 \to \B \big(L^2(K)_+ \big)$ by
\begin{eqnarray*}
F_a^1(x)&:=&\F(a)\big(\P^{+,x} \big)~~~\forall x \in i[0, \iy), \\
F_a^1(x)&:=&\F(a)\big({\cal C}^x \big)~~~\forall x \in (0, 1)~~~\text{and} \\
F_a^1(1)&:=&K_{(1)} \circ \F(a) \big(\D_1^+ \big) \circ K_{(1)}^{-1} \circ p_+~+~K_{[1]} \circ \F(a) \big(\D_1^- \big) \circ K_{[1]}^{-1} \circ p_-~+~ \F(a) \big(\F_1 \big) \circ p_0.
\end{eqnarray*}
By Property 1. of Theorem \ref{theo1 - sl2r}, $F_a^1(x) \in \K \big(L^2(K)_+ \big)$ for all $x \in I_1 \setminus \{1\}$. Moreover, since $\F(a) \big(\D_1^+ \big)$, $\F(a) \big(\D_1^- \big)$ and $\F(a) \big(\F_1 \big)$ are also compact, their composition with the bounded ope\-ra\-tors $K_{(1)}$, $K_{(1)}^{-1}$, $K_{[1]}$, $K_{[1]}^{-1}$, $p_+$, $p_-$ and $p_0$ is compact as well. Therefore, $F_a^1(1)\in \K \big(L^2(K)_+ \big)$ too. \\
By Property 4. of Theorem \ref{theo1 - sl2r}, $F_a^1$ vanishes at infinity. Moreover, for all $x \in I_1 \setminus \{0,1\}$, $F_a^1$ is obviously continuous in $x$. \\
For the continuity in $0$, let $\overline{u}=(u_j)_{j \in \N}$ be a sequence in $(0,1)$ converging to $0$. Then,
\begin{eqnarray*}
\lim \limits_{j\to\iy } \big\| F_a^1(u_j)-F_a^1(0)\big\|_{op}
&=&\lim \limits_{j\to\iy } \big\| \F(a) \big({\cal C}^{u_j} \big)-\F(a)\big( \P^{+,0} \big) \big\|_{op} \\
&=&\lim \limits_{j\to\iy } \big\| {\cal C}^{u_j}(a)- \P^{+,0}(a) \big\|_{op} \\
&=&\lim \limits_{j\to\iy } \big\| K_{u_j} \circ \P^{+,u_j}(a) \circ K_{u_j}^{-1}- \id \circ \P^{+,0}(a) \circ \id^{-1} \big\|_{op}~=~0, 
\end{eqnarray*}
as $K_{u_j} \overset{j \to \iy}{\longrightarrow} \id$ and with the same arguments as in the proof of (\ref{3afuergamma8}) in Section \ref{cond 3a}. \\
For the continuity in $1$, let $\overline{u}=(u_j)_{j \in \N}$ a sequence in $(0,1)$ converging to $1$. By Property 5. of Theorem \ref{theo1 - sl2r} and since $F_a^1(1)= \tilde{\nu}_{\overline{u},j} \big( \F(a) \big)= \tilde{\nu} \big( \F(a) \big)$ by the definition of $\tilde{\nu}$ in Case (iii) of Section \ref{cond 3b},
\begin{eqnarray*}
\lim \limits_{j\to\iy } \big\| F_a^1(u_j)-F_a^1(1)\big\|_{op}
~=~\lim \limits_{j\to\iy } \big\| \F(a) \big({\cal C}^{u_j} \big)-\tilde{\nu} \big(\F(a)\big)\big\|_{op}~=~0. 
\end{eqnarray*}
Hence, $F_a^1$ is also continuous in $0$ and in $1$ and thus $F_a^1 \in C_{\iy} \Big(i [0, \iy) \cup [0,1], \K \big(L^2(K)_+ \big)\Big)$. \\
Since $p_+$, $p_-$ and $p_0$ equal the identity on the image of $K_{(1)}$, $K_{[1]}$ and respectively $\F(a) \big(\F_1 \big)$, as discovered in Section \ref{cond 3b}, and since $p_+ \circ p_-=p_- \circ p_+=p_+ \circ p_0=p_0 \circ p_+=p_- \circ p_0=p_0 \circ p_-=0$, $F_a^1(1)$ commutes with $p_+$, $p_-$ and $p_0$. \\
~\\
On the other hand, taking a function $F_1 \in C_{\iy} \Big(i [0, \iy) \cup [0,1], \K \big(L^2(K)_+ \big)\Big)$ that commutes with $p_+$, $p_-$ and $p_0$, it has to be shown that there is an operator field $\va_{F_1}$ over $\wh{G}_{even}$ that meets the Properties 1. to 5. of Theorem \ref{theo1 - sl2r}: \\
Define
\begin{eqnarray*}
\va_{F_1}\big( \P^{+,x} \big)&:=&F_1(x) \in \B \big(L^2(K)_+ \big)~~~\forall x \in i[0, \iy), \\
\va_{F_1}\big({\cal C}^x \big)&:=&F_1(x) \in \B \big(L^2(K)_+ \big)~~~\forall x \in (0, 1), \\
\va_{F_1}\big( \D_1^{+} \big)&:=&K_{(1)}^{-1} \circ p_+ \circ F_1(1)  \circ K_{(1)} \in \B \big(\H_{(1)} \big), \\
\va_{F_1}\big( \D_1^{-} \big)&:=&K_{[1]}^{-1} \circ p_- \circ F_1(1)  \circ K_{[1]} \in \B \big(\H_{[1]} \big)~~~\text{and} \\
\va_{F_1} \big(\F_1 \big)&:=& p_0 \circ F_1(1) \circ p_0 \in \C.
\end{eqnarray*}
By the definition of $C_{\iy} \Big(i [0, \iy) \cup [0,1], \K \big(L^2(K)_+ \big)\Big)$ and as the composition of the compact operator $F_1(1)$ with bounded operators is compact again, the Properties 1. to 4. are obviously fulfilled. \\
For Property 5., there are two cases to consider: a sequence  in $(0,1)$ converging to $0$ and a sequence in $(0,1)$ converging to $1$. \\
So, first let $(u_j)_{j \in \N}$ be a sequence in $(0,1)$ converging to $0$. Then, by the definition of $\tilde{\nu}$ in Case (ii) of Section \ref{cond 3b}
\begin{eqnarray*}
\big\| \va_{F_1} \big({\cal C}^{u_j} \big)- \tilde{\nu} \big(\va_{F_1} \big) \big\|_{op}
&=& \big\| F_1(u_j) - \va_{F_1} \big(\P^{+,0} \big)   \big\|_{op}
~=~ \big\| F_1(u_j) -  F_1(0) \big\|_{op} \overset{j \to \iy}{\longrightarrow} 0,
\end{eqnarray*}
since $F_1$ is continuous in $0$. \\
Now, let $(u_j)_{j \in \N}$ a sequence in $(0,1)$ converging to $1$. Then, as $F_1(1)$ commutes with $p_+$, $p_-$ and $p_0$, by the definition of $\tilde{\nu}$ in Case (iii) of Section \ref{cond 3b} and since $p_+ + p_-+p_0= \id_{L^2(K)_+ \to L^2(K)_+}$,
\begin{eqnarray*}
&&\big\| \va_{F_1} \big({\cal C}^{u_j} \big)- \tilde{\nu} \big(\va_{F_1} \big) \big\|_{op}\\
&=& \Big\| F_1(u_j) - K_{(1)} \circ \va_{F_1} \big(\D_1^+ \big) \circ K_{(1)}^{-1} \circ p_+ -K_{[1]} \circ \va_{F_1} \big(\D_1^- \big) \circ K_{[1]}^{-1} \circ p_- - \va_{F_1} \big(\F_1 \big) \circ p_0 \Big\|_{op} \\
&=& \Big\| F_1(u_j) - K_{(1)} \circ K_{(1)}^{-1} \circ p_+ \circ F_1(1)  \circ K_{(1)} \circ K_{(1)}^{-1} \circ p_+\\
&&~- K_{[1]} \circ K_{[1]}^{-1} \circ p_- \circ F_1(1)  \circ K_{[1]} \circ K_{[1]}^{-1} \circ p_- - p_0 \circ F_1(1) \circ p_0 \circ p_0 \Big\|_{op} \\
&=& \big\| F_1(u_j) - p_+ \circ F_1(1) \circ p_+ -  p_- \circ F_1(1) \circ p_- -p_0 \circ F_1(1) \circ p_0 \big\|_{op} \\
&=& \big\| F_1(u_j) - F_1(1) \circ \big(p_+ +  p_- + p_0 \big) \big\|_{op} \\
&=& \| F_1(u_j) - F_1(1) \|_{op} \overset{j \to \iy}{\longrightarrow} 0
\end{eqnarray*}
because of the continuity of $F_1$ in $1$. \\
Therefore, Property 5 is also fulfilled. \\
~\\
Furthermore, one defines for all $a \in C^*(G)$ the function $F_a^2:I_2 \to \B \big(L^2(K)_- \big)$ by
\begin{eqnarray*}
F_a^2(x)&:=&\F(a)\big(\P^{-,x} \big)~~~\forall x \in i(0, \iy)~~~\text{and} \\
F_a^2(0)&:=&\F(a)\big(\D_+ \big) \circ p_+~+~\F(a)\big(\D_- \big) \circ p_-.
\end{eqnarray*}
Using the same arguments as above, one gets the desired properties of the function $F_a^2$ as well. \\

Next, take a function $F_2 \in C_{\iy} \Big(i [0, \iy), \K \big(L^2(K)_- \big)\Big)$ that commutes with $p_+$ and $p_-$. An operator field $\va_{F_2}$ over $\wh{G}_{odd}$ meeting the Properties 1. to 5. of Theorem \ref{theo1 - sl2r} has to be constructed: \\
Define
\begin{eqnarray*}
\va_{F_2} \big(\P^{-,x} \big)&:=&F_2(x) \in \B \big(L^2(K)_- \big)~~~\forall x \in i (0, \iy), \\
\va_{F_2} \big(\D_+ \big)&:=&p_+ \circ F_2(0) \in \B \big(\H_{\D_+}\big) ~~~\text{and} \\
\va_{F_2} \big(\D_- \big)&:=& p_- \circ F_2(0) \in \B \big(\H_{\D_-} \big).
\end{eqnarray*}
Here again, the proof of Properties 1. to 5. is similar to the one above. \\

Now, take the infinite-dimensional and sepa\-rable Hilbert space $\H_{\D}$ for the representations $\D_m^+$ and $\D_m^-$ for $m>1$, fixed in Chapter \ref{rep}. Then, define for every $a \in C^*(G)$ the function $F_a^3:I_3 \to \B \big(\H_{\D} \big)$ by
\begin{eqnarray*}
F_a^3(x)&:=&\F(a)\big( \D_x^+ \big)~~~\forall x \in \Z_{>1}~~~\text{and} \\
F_a^3(x)&:=&\F(a)\big( \D_{-x}^- \big)~~~\forall x \in \Z_{<-1}.
\end{eqnarray*}
Here, Property 5. of Theorem \ref{theo1 - sl2r} does not emerge and the Properties 1. to 4. are obvious.\\
Taking a function $F_3 \in C_{\iy}\Big(\Z \setminus \{-1,0,1\}, \K \big(\H_{\D} \big) \Big)$, one has to choose
\begin{eqnarray*} 
\va_{F_3} \big( \D_x^+ \big)&:=&F_3(x) \in \B \big(\H_{\D} \big)~~~\forall x \in \Z_{>1}~~~\text{and} \\
\va_{F_3} \big( \D_x^- \big)&:=&F_3(-x) \in \B \big(\H_{\D} \big)~~~\forall x \in \Z_{>1}
\end{eqnarray*}
and it is again easy to check that $\va_F$ complies with the properties of Theorem \ref{theo1 - sl2r}. \\
\qed

\section{Acknowledgements} 
This work is supported by the Fonds National de la Recherche, Luxembourg (Project Code 3964572). 

\section{References}
\begin{enumerate}
\item \label{bel-bel-lud} I.Beltita, D.Beltita and J.Ludwig, Fourier transforms of $C^*$-algebras of nilpotent Lie groups, arXiv:1411.3254, 2014.
\item \label{cohn} L.Cohn, Analytic Theory of the Harish-Chandra C-Function, Springer-Verlag, Berlin-Heidelberg-New York, 1974.
\item \label{dix} J.Dixmier, $C^*$-algebras. Translated from French by Francis Jellett. North-Holland Mathe\-matical Library, Vol. 15, North-Holland Publishing Co., Amsterdam-New York-Oxford, 1977.
\item \label{janne} J.-K.G\"unther and J.Ludwig, The $C^*$-algebras of connected real two-step nilpotent Lie groups, Revista Matem\'atica Complutense 29(1), pp. 13-57, 10.1007/s13163-015-0177-7, 2016.
\item \label{knapp} A.Knapp, Representation Theory of Semisimple Groups. An Overview Based on Examples. Princeton University Press, Princeton, New Jersey, 1986.
\item \label{lang} S.Lang, $SL_2(\R)$, New York, 1985.
\item \label{lin-ludwig} Y.-F.Lin and J.Ludwig, The $C^*$-algebras of $ax+b$-like groups, Journal of Functional Analy\-sis 259, pp. 104-130, 2010.
\item \label{ludwig-turowska} J.Ludwig and L.Turowska, The $C^*$-algebras of the Heisenberg Group and of thread-like Lie groups, Math. Z. 268, no. 3-4, pp. 897-930, 2011. 
\item \label{hedidr} H.Regeiba, Les $C^*$-alg\`{e}bres des groupes de Lie nilpotents de dimension $\leq 6$, Ph.D. thesis at the Universit\'{e} de Lorraine, 2014.
\item \label{hedi} H.Regeiba and J.Ludwig, $C^*$-Algebras with Norm Controlled Dual Limits and Nilpotent Lie Groups, arXiv:1309.6941, 2013.
\item \label{wal} N.Wallach, Real Reductive Groups I. Academic Press, Pure and Applied Mathematics, San Diego, 1988.
\item \label{walII} N.Wallach, Real Reductive Groups II. Academic Press, Pure and Applied Mathematics, San Diego, 1992.
\item \label{was} A.Wassermann, Une d\'emonstration de la conjecture de Connes-Kasparov pour les groupes de Lie lin\'eaires connexes r\'eductifs, C.R. Acad. Sci. Paris S\'erie I Math. 304, no.18, pp. 559-562, 1987 (French with English summary).
\end{enumerate} 

\end{document}